\documentclass[journal ]{new-aiaa}
\usepackage[utf8]{inputenc}

\usepackage{amsmath}
\usepackage{bm}
\usepackage{graphicx}
\usepackage[version=4]{mhchem}
\usepackage{siunitx}
\usepackage{longtable,tabularx}
\usepackage{gensymb}
\usepackage{mathtools}
\usepackage{booktabs}
\usepackage{float}
\usepackage[ruled, vlined, linesnumbered]{algorithm2e}
\usepackage{subcaption}
\usepackage[flushleft]{threeparttable}
\usepackage{enumitem}
\usepackage{multirow}
\usepackage{lscape}
\usepackage{adjustbox}
\usepackage{mdframed}
\usepackage{pdflscape}
\usepackage{color,soul}
\usepackage[color=gray!30]{todonotes}
\usepackage{colortbl}
\definecolor{myblue}{rgb}{0, 0.23, 0.64}

\usepackage{hyperref}
\hypersetup{
    colorlinks=true,
    linkcolor=myblue,
    filecolor=magenta,
    citecolor=myblue,
    urlcolor=cyan,
}
\usepackage{dsfont}
\usepackage[flushleft]{threeparttable}
\SetKwComment{Comment}{$\triangleright$\ }{}

\usepackage{amsthm}
\theoremstyle{definition}

\newtheorem{definition}{Definition}
\newtheorem{example}{Example}
\newtheorem{formulation}{Formulation}

\theoremstyle{plain}
\newtheorem{theorem}{Theorem}

\theoremstyle{remark}
\newtheorem{remark}{Remark}

\newcommand{\LR}{(\hyperlink{(LR)}{LR})\xspace}
\newcommand{\LRo}{(\hyperlink{(LR1)}{LR1})\xspace}
\newcommand{\LRt}{(\hyperlink{(LR2)}{LR2})\xspace}

\let\svthefootnote\thefootnote
\newcommand\freefootnote[1]{%
  \let\thefootnote\relax%
  \footnotetext{#1}%
  \let\thefootnote\svthefootnote%
}

\author{Hang Woon Lee\footnote{Assistant Professor, Department of Mechanical and Aerospace Engineering; hangwoon.lee@mail.wvu.edu. Member AIAA (Corresponding Author).}}
\affil{West Virginia University, Morgantown, WV, 26506}
\author{Koki Ho\footnote{Associate Professor, Daniel Guggenheim School of Aerospace Engineering, Senior Member AIAA.}}
\affil{Georgia Institute of Technology, Atlanta, GA, 30332}

\begin{document}

\freefootnote{This work was presented at the AAS/AIAA Astrodynamics Specialist Conference, Virtual, August 9-11, 2021 as Paper AAS 21-719.}

\title{Regional Constellation Reconfiguration Problem: Integer Linear Programming Formulation and Lagrangian Heuristic Method}
\maketitle

\begin{abstract}
A group of satellites, with either homogeneous or heterogeneous orbital characteristics and/or hardware specifications, can undertake a reconfiguration process due to variations in operations pertaining to Earth observation missions. This paper investigates the problem of optimizing a satellite constellation reconfiguration process against two competing mission objectives: (i) the maximization of the total coverage reward and (ii) the minimization of the total cost of the transfer. The decision variables for the reconfiguration process include the design of the new configuration and the assignment of satellites from one configuration to another. We present a novel bi-objective integer linear programming formulation that combines constellation design and transfer problems. The formulation lends itself to the use of generic mixed-integer linear programming (MILP) methods such as the branch-and-bound algorithm for the computation of provably-optimal solutions; however, these approaches become computationally prohibitive even for moderately-sized instances. In response to this challenge, this paper proposes a Lagrangian relaxation-based heuristic method that leverages the assignment problem structure embedded in the problem. The results from the computational experiments attest to the near-optimality of the Lagrangian heuristic solutions and a significant improvement in the computational runtime compared to a commercial MILP solver.
\end{abstract}

\section*{Nomenclature}
{\renewcommand\arraystretch{1.0}
\noindent\begin{longtable*}{@{}l @{\quad=\quad} l@{}}
    \multicolumn{2}{@{}l}{Orbital Elements}\\
    $a$ & semi-major axis \\
    $e$ & eccentricity \\
    $inc$ & inclination \\
    $M$ & mean anomaly \\
    $u$ & argument of latitude \\
    $\omega$ & argument of periapsis \\
    $\Omega$ & right ascension of the ascending node \\
    \multicolumn{2}{@{}l}{Parameters and Decision Variables}\\
    $b$ & coverage timeline \\
    $c$ & assignment cost \\
    $r$ & coverage threshold \\
    $T$ & mission planning horizon period \\
    $v$ & visibility profile \\
    $V$ & visibility matrix \\
    $x$ & constellation pattern variable \\
    $y$ & coverage state variable \\
    $Z$ & objective function value \\
    $\varphi$ & assignment variable \\
    $\pi$ & coverage reward \\
    $\varepsilon$ & epsilon-constraint method parameter \\
    $\lambda$ & Lagrange multiplier \\
    $\vartheta$ & elevation angle \\
    $\theta$ & subgradient method step size \\
    \multicolumn{2}{@{}l}{Sets and Indices}\\
    $\mathcal{E}$ & set of edges \\
    $\mathcal{G}$ & graph \\
    $\mathcal{I}$ & set of satellites (index $i$) \\
    $\mathcal{J}$ & set of orbital slots (index $j$) \\
    $\mathcal{N}$ & exchange neighborhood \\
    $\mathcal{P}$ & set of target points (index $p$) \\
    $\mathcal{S}$ & set of subconstellations (index $s$) \\
    $\mathcal{T}$ & set of time steps (index $t$) \\
    $\mathbb{R}_{\geq0}$ & set of non-negative real numbers \\
    $\mathbb{Z}_{\geq0}$ & set of non-negative integers \\
\end{longtable*}}

\section{Introduction}
Satellite constellation systems often face varying mission requirements and environments during their operations. These variations may arise from changes in the area of interest (e.g., disaster monitoring \cite{he2020reconfigurable}, temporary reconnaissance \cite{chen2015reconfiguration}, and theater situational awareness missions), or from modifications to the desired coverage performance, such as switching from sporadic coverage to uninterrupted, single-fold coverage. Moreover, the systems themselves may need to adapt due to the addition of new satellites, for example, through staged deployment \cite{deweck2004staged,arnas2022uniform}, or the loss of existing satellites due to failures \cite{ferringer2009many} and/or end-of-life decommissions. Under such circumstances, it is logical for system operators to consider options to ``reconfigure'' an existing constellation system to maximize the utility of active on-orbit assets rather than launching an entirely new constellation.

We define constellation reconfiguration as the process of transforming an existing configuration into another to maintain the system in an optimal state, given a set of new mission requirements \cite{davis2010constellation,fakoor2016}. The design of a reconfiguration process is nontrivial and involves interdisciplinary fields of studies, such as satellite constellation design theory, orbital transfer trajectory optimization, and mathematical programming, to enable a robust constellation reconfiguration framework.

Of particular interest to this paper is the topic of satellite constellation reconfiguration in the context of Earth observations (EO). Many present-day EO satellite systems are monolithic or small-scale constellation systems distributed in near-polar low Earth orbits (LEO), mostly in sun-synchronous orbits to leverage consistent illumination conditions. Near-polar orbits enable EO satellite systems to scan different parts of the globe in each orbit, making them ideal for detecting changes in the Earth's land cover, vegetation, and civil infrastructures. However, the long revisit time for a particular target makes near-polar orbits unsuitable for missions that require rapid adaptive mission planning and enhanced coverage, such as satellite-based emergency mapping, surveillance, and reconnaissance missions, to name a few \cite{denis2016,voigt2016global}. Of latest attention to the EO community has been the concept of agile satellites that assume attitude control capability, which is deemed to enhance the overall system responsiveness and scheduling efficiency \cite{wang2021agile}. Recently, several works have explored the concept of maneuverable satellites in the domain of EO satellite systems as a new paradigm to bolster the system observation capacity by directly manipulating the orbits \cite{paek2019optimization,mcgrath2019general,he2020reconfigurable,zhang2021,morgan2023}. In this paper, we investigate the concept of reconfiguration as a means for system adaptability and responsiveness that adds a new dimension to the operation of next-generation EO satellite constellation systems.

The problem of satellite constellation reconfiguration consists of two different, yet coupled, problems: the \textit{constellation design problem} and the \textit{constellation transfer problem} \cite{davis2010constellation,appel2014optimization,legge2014optimization}. The former deals with the optimal design of a (destination) constellation configuration that satisfies a set of mission requirements; the latter is concerned with the minimum-cost transportation of satellites from one configuration to another provided the knowledge of both end states. Although we may approach these two interdependent problems independently in a sequential manner (i.e., a destination configuration is first designed and followed by the optimal assignment of satellites to new orbital slots), the outcome of such an open-loop procedure may result in a suboptimal reconfiguration process as a whole \cite{davis2010constellation,legge2014optimization}. Without taking into account the satellite transportation aspect in design, the optimized new configuration may be too costly or, in fact, infeasible to achieve. This background motivates us to concurrently consider constellation design and transfer aspects in satellite constellation reconfiguration.

The problem of concurrent constellation design and transfer optimization is highly complex and challenging. While it is well known that the constellation transfer problem can be formulated as an assignment problem \cite{deweck2008optimal}, formulating the constellation design problem faces a unique mathematical programming challenge due to (i) the potentially complex (spatiotemporally-varying) regional coverage and (ii) the reconfiguration problem with the cardinality constraint. First, the constellation design problem for complex regional coverage may need to incorporate design attributes such as heterogeneity among member satellites (e.g., different orbits and hardware specs) and asymmetry in satellite distributions. The classical constellation patterns, such as the streets-of-coverage \cite{luders1961,luders1974}, Walker patterns \cite{walker1970,walker1977,walker1984}, and the tetrahedron elliptical constellation \cite{draim1987common}, are limiting due to symmetry and sparsity in satellite distribution, especially with a small number of satellites. Ref.~\cite{lee2020satellite} has shown that, for complex regional coverage, relaxing symmetry and homogeneity assumptions of the classical methods enables the exploration of larger design space and hence leads to the discovery of more efficient constellation pattern sets. However, incorporating asymmetric patterns into the constellation reconfiguration while considering the transfer cost remains a challenging problem, unaddressed in the literature. Second, in the context of satellite constellation reconfiguration, the design of a destination configuration may be restricted to a given number of satellites, which we refer to as the cardinality constraint. This is a logical assumption to make because, without the enforcement of the cardinality constraint, the optimal design of a destination configuration may require substantially more satellites than are readily available for orbital maneuvers. Launching a set of new satellites to fulfill this deficit within a limited time window can be challenging from a financial and operational perspective. Therefore, addressing the reconfiguration problem while adhering to the cardinality constraint adds an extra layer of complexity.

The challenge in solving the satellite constellation reconfiguration problem lies not only in integrating design and transfer aspects but also in devising a solution approach that is computationally efficient and yields high-quality solutions, particularly for mission scenarios that require a rapid system response. The way we formulate the integrated design-transfer model determines its mathematical properties and the pool of applicable solution algorithms, which in turn affects the time complexity of retrieving solutions. Several satellite constellation reconfiguration studies have been conducted, covering both constellation design and transfer in various problem settings \cite{ferringer2009many,chen2015reconfiguration,paek2019optimization,he2020reconfigurable,zhu2010satellite}. These studies have demonstrated the value of concurrent optimization, but their mathematical problem formulations are generally nonlinear and often employ meta-heuristic algorithms. While these algorithms can be efficient in obtaining high-quality solutions, they can be computationally expensive for highly-constrained problems and cannot certify the optimality (or the optimality gap) of the obtained solutions. Therefore, the principal challenge we face in this work involves streamlining the entire pipeline, from formulating the constellation design problem to integrating constellation design and transfer aspects and developing a computationally-efficient solution method.

The contributions of this paper are as follows. We present an integer linear program (ILP) formulation of the design-transfer problem, referred to as the \textit{Regional Constellation Reconfiguration Problem} (RCRP). This formulation incorporates both constellation design and constellation transfer aspects, which are typically considered independent and serial in current state-of-the-art techniques. The RCRP utilizes the maximal covering location problem formulation found in facility location problems for constellation design, and the assignment problem for constellation transfer, both of which are ILPs. By integrating these two problems, a larger design space is explored and operators are provided with trade-off analysis between transportation cost and coverage performance. The proposed model supports various mission concepts of operations that arise in regional coverage missions. The presented RCRP formulation enables the use of using mixed-integer linear programming (MILP) methods, such as the branch-and-bound algorithm, to obtain globally-optimal reconfiguration solutions. However, this approach becomes intractable for moderately-sized instances. To address this challenge, a Lagrangian relaxation-based solution method is proposed to approach large-scale optimization. This method relaxes a set of constraints to reveal and exploit the special substructure of the problem, making it easier to solve. The results of the computational experiments demonstrate the near-optimality of the Lagrangian heuristic solutions, compared to solutions obtained by a commercial solver, with significantly faster runtime.

The remainder of this paper is organized as follows. Section~\ref{sec:background} provides an overview of the constellation-coverage model and discusses the optimization formulations for constellation design and transfer problems. Section~\ref{sec:integratedmodel} presents a mathematical formulation of the integrated design-transfer problem and examines its characteristics. Section~\ref{sec:solution} introduces the developed Lagrangian relaxation-based solution method for addressing the proposed problem formulation. Section~\ref{sec:experiments} conducts computational experiments to demonstrate the effectiveness of the developed method and provides an illustrative example applied to the case of federated disaster monitoring. Finally, Section~\ref{sec:conclusion} concludes this paper.

\section{Constellation Design and Transfer Problems} \label{sec:background}

In this section, we construct a constellation-coverage model (Section~\ref{sec:model}), propose an optimization problem formulation for the constellation design problem (Section~\ref{sec:mcp}), and review a constellation transfer problem model by Ref.~\cite{deweck2008optimal} (Section~\ref{sec:ap}). The materials discussed in this section lay the foundation for the integrated design-transfer model in Section~\ref{sec:integratedmodel}.

Several comments on the notations are noted. The asterisk symbol in superscript $(\cdot)^\ast$ denotes the optimality of a variable $(\cdot)$. $Z(\cdot)$ denotes the optimal objective function value of a given problem with parameters $(\cdot)$. $Z_{\text{LP}}$ denotes the optimal value of a given problem with integrality constraints dropped, hence the name linear programming (LP) relaxation bound. $\text{Co}(\cdot)$ denotes the convex hull of a set $(\cdot)$, and $|\cdot|$ denotes the cardinality of a set $(\cdot)$.

\subsection{Constellation-Coverage Model} \label{sec:model}
We introduce the constellation-coverage model that relates the configuration of a constellation system with its coverage performance. In this model, the finite time horizon of period $T$ is discretized into a set of time steps with a step size $\Delta t$. Let $\mathcal{T} \coloneqq \{0,1,\dots,m-1\}$, where $m\Delta t= T$, be the set of time step indices $t$ such that the set $\{t(\Delta t):t\in\mathcal{T}\}$ is the discrete-time finite horizon. The set of orbital slot indices is denoted by $\mathcal{J}$. Each orbital slot $j\in\mathcal{J}$ is defined by a unique set of orbital elements $\textbf{\oe}_j=(a_j,e_j,inc_j,\omega_j,\Omega_j,M_j)$. Here, $a$, $e$, $inc$, $\omega$, $\Omega$, and $M$ each represents the semi-major axis, eccentricity, inclination, argument of periapsis, right ascension of ascending node (RAAN), and mean anomaly of an orbit, respectively. For circular orbits, we use the argument of latitude $u$. We also let $\mathcal{P}$ be the set of target point indices $p$.

\subsubsection{Model Definitions}
For ease of description, and without loss of generality, we consider the model for a single target point.

\begin{definition}[Visibility matrix]
Let $V_{tj}$ denote the Boolean visibility state, which equals 1 if a satellite in orbital slot $j$ covers the target point at time step $t$. We let $\bm{V}=(V_{tj}\in\mathbb{Z}_{\geq0}^2:t\in\mathcal{T},j\in\mathcal{J})$ denote a visibility matrix where $\mathbb{Z}_{\geq0}$ denotes the set of non-negative integers.
\end{definition}

To construct $\bm{V}$, the following parameters need to be specified for each orbital slot: the orbital elements $\textbf{\oe}_j$ at the epoch, the minimum elevation angle threshold $\vartheta_{\min}$ for a target point (and/or the field-of-view of a satellite sensor), the coordinates of a target point, and the epoch at which the finite time horizon is referenced. With these parameters, the orbital slot is numerically propagated under the governing equations of motion (e.g., $J_2$-perturbed two-body motion) for a finite time horizon of period $T$. At each time step $t$, a Boolean visibility masking is applied to construct an element of the visibility matrix, $V_{tj}$.

\begin{definition}[Constellation pattern vector] \label{def:x}
A constellation pattern vector $\bm{x}=(x_j \in \{0,1\}:j\in\mathcal{J})$ specifies the relative distribution of satellites in a given system (or simply, the configuration of a constellation system). Each element of $\bm{x}$ is defined as:
\begin{equation*}
x_j\coloneqq\begin{cases}
1, &\text{if a satellite occupies orbital slot } j \\
0, &\text{otherwise}
\end{cases}
\end{equation*}
\end{definition}

\begin{definition}[Coverage timeline]
Let $b_t$ be the number of satellite(s) in view from the target point at time step $t$. Then, we let $\bm{b}=(b_t\in\mathbb{Z}_{\geq0}:t\in\mathcal{T})$ denote a coverage timeline. Here, the visibility of a satellite from a target point follows from the Boolean visibility masking.
\end{definition}

\begin{remark}[Linear property]
We can relate visibility matrix $\bm{V}$, constellation pattern vector $\bm{x}$, and coverage timeline $\bm{b}$ as a linear system. Mathematically,
\begin{equation}
\label{eq:Vxb}
b_t=\sum_{j\in \mathcal{J}} V_{tj} x_j
\end{equation}
\end{remark}

In this model, the set $\mathcal{J}$ can comprise orbital slots with different orbital characteristics without any predefined rule. Because satellites in these orbital slots experience different degrees of orbital perturbations over time, the constellation-coverage model is only valid within the specified time horizon of period $T$. There will be a loss of fidelity in the constellation and coverage relationship beyond the specified time horizon. Such a case is indeed suitable for the planning of many temporary mission operations. However, some cases require persistent coverage of a region of interest for a long-term horizon. To account for it, we make several assumptions about the constellation-coverage model. In what follows, we review the assumptions and the definitions of the Access-Pattern-Coverage (APC) decomposition model by Ref.~\cite{lee2020satellite}.

\subsubsection{Special Case: APC Decomposition}
To guarantee persistent regional coverage, Ref.~\cite{lee2020satellite} introduced a particular constellation-coverage model called the APC decomposition (named after the three finite discrete-time sequences of the special case model: the visibility (Access) profile, constellation Pattern vector, and Coverage timeline) by making two assumptions about the set of orbital slots $\mathcal{J}$: (i) the repeating ground track (RGT) orbits and (ii) the common ground track constellation. Specifically, the conditions are:
\begin{enumerate}[label=(\roman*)]
    \item A ground track is the trace of a satellite's sub-satellite points on the surface of a planetary body. A satellite on an RGT orbit makes $N_\text{P}$ number of revolutions in $N_\text{D}$ number of nodal periods. There is a finite time horizon of period $T$ (often called a period of repetition) during which a satellite repeats its closed relative trajectory exactly and periodically. Expressing this condition \cite{mortari2004}, we get:
\begin{equation*}
T=N_\text{P}T_\text{S}=N_\text{D}T_\text{G}
\end{equation*}
where $N_\text{P}$ and $N_\text{D}$ are positive integers. $T_\text{S}$ is the nodal period of a satellite due to both nominal motion and perturbations and $T_\text{G}$ is the nodal period of Greenwich.
\item All satellites in a common ground track constellation share identical semi-major axis $a$, eccentricity $e$, inclination $inc$, and argument of periapsis $\omega$ but each satellite $i$ independently holds a pair of right ascension of ascending node (RAAN) $\Omega_i$ and initial mean anomaly $M_i$ that satisfy the following distribution rule \cite{avendano2013}:
\begin{equation}
\label{eq:raan_m}
N_{\text{P}} \Omega_i+N_{\text{D}} M_i = \text{constant} \ \text{mod} \ 2\pi
\end{equation}
\end{enumerate}

With these assumptions, we add the following new definitions to the model to accommodate the special case.

\begin{definition}[Reference visibility profile]
Let $v_t$ denote the Boolean visibility state that equals 1 if a reference satellite covers a target point at time step $t$ (0 otherwise). Then, we denote $\bm{v}=(v_{t} \in \{0,1\}:t\in\mathcal{T})$ the reference visibility profile.
\end{definition}

\begin{definition}[Visibility circulant matrix] \label{def:circulant}
A visibility circulant matrix $\bm{V}$ is the $m\times m$ matrix whose columns are the cyclic permutations of $\bm{v}$:
\begin{equation*}
\renewcommand{\arraystretch}{0.7}
    \bm{V}= \text{circ}(\bm{v}) =
    \begin{bmatrix}
    v_0 & v_{m-1} & \cdots & v_1 \\
    v_1 & v_0 & \cdots & v_2 \\
    \vdots & \vdots & \ddots & \vdots \\
    v_{m-1} & v_{m-2} & \cdots & v_0
    \end{bmatrix}
\end{equation*}
where the $(t,j)$ entry of $\bm{V}$ is denoted with the modulo operator as $V_{tj}=v_{(t-j) \bmod m}$; $\text{circ}(\cdot)$ is the circulant operator that takes $\bm{v}$ as the argument and generates a circulant matrix as defined above.
\end{definition}

\begin{remark}[Circular convolution operation \cite{lee2020satellite}]
We can relate reference visibility profile $\bm{v}$, constellation pattern vector $\bm{x}$, and coverage timeline $\bm{b}$ in the manner prescribed by a \textit{circular convolution operation}. Mathematically,
\begin{equation}
\label{eq:circular}
b_t=\sum_{j\in \mathcal{J}} v_{(t-j) \bmod m} x_j
\end{equation}
\end{remark}

Following from the definition of the ($t,j$) entry of $\bm{V}$ in Definition~\ref{def:circulant}, Eq.~\eqref{eq:circular} can be written as a linear system in terms of a reference visibility circulant matrix: $\bm{b}=\bm{V}\bm{x}$, which is in the form of Eq.~\eqref{eq:Vxb}. 

There are two notable benefits to this special case. One unique advantage of this special case is that it only requires knowledge of the reference visibility profile and the distribution of satellites along the common relative trajectory to quantify the satellite coverage state $b_t$ of a target point at time step $t$. The construction of $\bm{V}$ is significantly faster than the generic case due to the use of the circulant operator. Another advantage is that, as will be discussed later in this paper, having $\bm{V}$ as a circulant matrix can lead to useful mathematical properties. In particular, in Section~\ref{sec:mcp}, we show that the upper bound of the LP relaxation of the maximum coverage problem can be analytically computed if $\bm{V}$ is circulant. Although not directly relevant to the main contents of this paper, circulant matrices can be used to leverage efficient solution methods for certain classes of problems (e.g., set covering problems with circulant matrices, as demonstrated in references \cite{bartholdi1980,bartholdi1981}).

\subsection{Constellation Design: Maximum Coverage Problem} \label{sec:mcp}
In this subsection, we introduce the Maximum Coverage Problem (MCP), which models the constellation design aspect of a constellation reconfiguration process. The MCP is based on the constellation-coverage model presented in Section~\ref{sec:model} and serves as one of two essential components of the proposed RCRP formulation, which will be introduced in Section~\ref{sec:integratedmodel}. Additionally, we elucidate some of the interesting properties of this MCP formulation.

Consider a problem setting where $\mathcal{T}$ is the set of time step indices and $\mathcal{J}$ is the set of orbital slot indices. Without loss of generality, we consider the problem for a single target point of interest. Given a finite time horizon of period $T$, the time-dependent observation reward $\bm{\pi}=(\pi_t\in\mathbb{R}_{\geq0}:t\in\mathcal{T})$ is defined for a target point of interest. The goal is to locate $n$ satellites in $\mathcal{J}$ such that the total observation reward obtained by covering the target point is maximized.

To obtain the reward $\pi_t$ at time step $t$, the target point must be covered. The target point is considered covered if there is at least $r_t$ satellite(s) in view; the positive time-dependent coverage threshold $\bm{r}=(r_t\in\mathbb{Z}_{\geq0}:t\in\mathcal{T})$ is a user-supplied parameter vector. We can model this system by defining two sets of decision variables: the constellation pattern variables $\bm{x}$ and coverage state variables $\bm{y}=(y_t\in\{0,1\}:t\in\mathcal{T})$, and a set of inequalities that links $\bm{x}$ and $\bm{y}$. The decision variable $x_j=1$ if a satellite occupies orbital slot $j$ ($x_j=0$ otherwise; see Definition~\ref{def:x}). Each element $y_t$ of the coverage state variables takes the value of unity if and only if the coverage threshold of the target point is satisfied at time step $t$ ($y_t=0$ otherwise). Mathematically,
\begin{equation}
\label{eq:mcpy}
    y_t = \begin{cases}
1, &\text{if } b_t=\sum_{j\in \mathcal{J}}V_{tj}x_j\ge r_t \\
0, &\text{otherwise}
\end{cases}
\end{equation}

As can be seen in Eq.~\eqref{eq:mcpy}, the coverage state of the target point is conditionally dependent on the configuration of a constellation system. We can linearize this relationship by introducing a set of inequalities that links $\bm{x}$ and $\bm{y}$ for all $t\in\mathcal{T}$:
\begin{equation*}
\sum_{j\in\mathcal{J}} V_{tj}x_j\ge r_t y_t, \quad \forall t\in \mathcal{T}
\end{equation*}

With these conditions as constraints, MCP aims to maximize the coverage reward of a satellite configuration with a given number of satellites, $n$. The preliminary version of the mathematical formulation of MCP was introduced in our earlier work \cite{Lee2020binary}, which employs the big-M method to linearize the conditional constraints. We present an improved formulation of MCP that achieves a tighter integrality gap.

\begin{formulation}[Maximum coverage problem] \label{md:mcp}
MCP is formulated as an integer linear program:
\begin{subequations}
\begin{alignat}{2}
	\text{(MCP)} \quad Z = \max
	\quad & \sum_{t\in \mathcal{T}}{\pi_ty_t} \label{eq:mcpo} \\
	\text{s.t.} \quad
	& \sum_{j\in\mathcal{J}} V_{tj}x_j\ge r_t y_t, \quad & \forall t\in \mathcal{T} \label{eq:mcpa} \\
	& \sum_{j\in \mathcal{J}} x_j = n \label{eq:mcpb} \\
	& x_j \in \{0,1\}, & \forall j\in \mathcal{J} \label{eq:mcpc} \\
	& y_t \in \{0,1\}, \label{eq:mcpd} & \forall t\in \mathcal{T} 
\end{alignat}
\end{subequations}
where $Z$ denotes the optimal value of MCP.
\end{formulation}

The objective function~\eqref{eq:mcpo} maximizes the total reward earned by covering a given target point. Constraints~\eqref{eq:mcpa} couples the configuration of a system to its coverage state of the target point. Constraint~\eqref{eq:mcpb} is the cardinality constraint that restricts the number of satellites to a fixed value of $n$. Constraints~\eqref{eq:mcpc} and \eqref{eq:mcpd} define the domains of decision variables.

The 0-1 integrality constraint on $y_t$ can be relaxed to $0\le y_t \le 1$ when $r_t=1$. This avoids the use of unnecessary integer definitions and could potentially facilitate the branch-and-bound algorithm necessary to find the optimal solution \cite{revelle2008}. Note that for $r_t >1$, $y_t$ may take a fractional value (for instance, if $r_t=2$ and $\sum_{j\in \mathcal{J}}V_{tj}x_j=1$, then $y_t$ can take 0.5), therefore the integrality constraints on $y_t$ must be enforced.

Note that for $r_t=1,\forall t \in \mathcal{T}$, the MCP can be shown equivalent to the maximal covering location problem (MCLP) that emerges in many problem contexts such as the facility location problem. The MCLP seeks to locate a number of facilities such that the weighted coverage of demand nodes is maximized; each facility is pre-specified with a service radius to which it can provide coverage. We can use the following analogy: the satellites are the facilities and the time steps are the demand nodes. (Conversely, this suggests the general varying-radius $\bm{r}$-fold coverage formulation of MCLP.) Unfortunately, the equivalence in the formulations informs us that MCP is NP-hard because of the NP-hardness of the MCLP \cite{megiddo1983}, which can be deduced using the argument of the reduction from the MCLP to the MCP. For more information on the mathematical formulation and the applications of the MCLP, readers are encouraged to refer to the original study by Church and ReVelle \cite{church1974}.

Expressing $Z_{\text{LP}}$ in terms of an optimal LP solution ($\bm{x}^\ast,\bm{y}^\ast$), we get:
\begin{equation}
    Z_{\text{LP}}=\sum_{t\in \mathcal{T}} \pi_t y_t^\ast \label{eq:zlp}
\end{equation}
where $y_t^\ast$ is determined from $x_j^\ast$ as
\begin{equation}
\label{eq:yt}
    y_t^\ast=\min\Bigg(\frac{1}{r_t}\sum_{j\in \mathcal{J}}V_{tj}x_j^\ast,1\Bigg)
\end{equation}
Equation~\eqref{eq:yt} follows from the fact that MCP is a maximization problem---$y_t$ variables will take their maximum values as bounded by Constraints~\eqref{eq:mcpa}. The second argument in the $\min(\cdot)$ operator bounds the maximum of $y_t$ to one, conforming with Eq.~\eqref{eq:mcpy}.

The discussion on the LP relaxation bound of MCP will be revisited later in Section~\ref{sec:solution} for the proposed solution method. Therefore, we provide additional implications for $Z_{\text{LP}}$ in this subsection to ensure completeness. First, we notice that $y_t^\ast$ in Eq.~\eqref{eq:zlp} requires knowledge of the optimal LP solution $x_j^\ast$. However, under special conditions, we can closely approximate $Z_{\text{LP}}$ by computing the upper bound $\hat{Z}_{\text{LP}}$ such that no knowledge of $x_j^\ast$ is needed. To show this, we derive $\hat{Z}_{\text{LP}}$ from Eq.~\eqref{eq:zlp} by moving the summation inside the $\min(\cdot)$ function. Expanding the first argument further to separate the first column of $\bm{V}$ and defining $\xi_t \coloneqq \pi_t/r_t$ as the reward-to-requirement ratio, we obtain Eq.~\eqref{eq:zlpb}:
\begin{align}
    \hat{Z}_{\text{LP}} &\coloneqq \min\Bigg(\underbrace{n\sum_{t\in\mathcal{T}} \xi_t V_{t1}}_{(1)}+\underbrace{\sum_{j\in\mathcal{J}\setminus \{1\}} x_j^\ast\Bigg[\sum_{t\in\mathcal{T}} \xi_t V_{tj}-\sum_{t\in\mathcal{T}} \xi_t V_{t1}\Bigg]}_{(2)},\sum_{t\in\mathcal{T}} \pi_t \Bigg) \label{eq:zlpb} \\
    &\ge Z_{\text{LP}} \nonumber
\end{align}
where we group the first argument of the $\min(\cdot)$ function into two terms.

If $\xi_t = \xi$ for all $t \in \mathcal{T}$ and $\bm{V}$ is a circulant matrix, then the terms within the bracket in Term~(2) cancel out. Most problem instances we deal with in this paper assume $r_t=r, \forall t \in \mathcal{T}$ (time-invariant, $r$-fold continuous coverage) and $\pi_t=\pi, \forall t \in \mathcal{T}$ (uniform coverage reward) such that Term~(2) vanishes. Therefore, we can conveniently express $\hat{Z}_{\text{LP}}=\min(n\sum_{t\in\mathcal{T}} \xi_t v_{t},\sum_{t\in\mathcal{T}}\pi_t)$ as a function of known parameters $n$, $\xi$, and $\bm{v}$ and without the needing to run LP.

\begin{example}[5-satellite MCP] \label{eg:3}
Let $\textbf{\oe}_0=(a,e,inc,\Omega,u)=(\SI{12758.5}{km},0,50\degree,50\degree,0\degree)$ be the orbital elements of the reference satellite defined in the J2000 frame. This corresponds to the RGT ratio of $N_\text{P}/N_\text{D}=6/1$, that is, a satellite makes six revolutions in one nodal day. Assume a single target point of interest $p$ with the geodetic coordinate $(40\degree \text{N},100\degree \text{W})$. The minimum elevation angle threshold $\vartheta_{\min}$ for the target is set to \SI{10}{deg}. Suppose we wish to maximize the coverage over this target point with five satellites. For the MCP-specific parameters, we let $r_t=1,\forall t \in \mathcal{T}$ and $\pi_t=1,\forall t \in \mathcal{T}$; this simplifies the MCP to the coverage percentage maximization problem. Solving the MCP to optimality, we get the optimum of $Z=398$, which translates into 79.6\% temporal coverage of target point $p$ by the optimal five-satellite configuration during the given repeat period $T$. Interpreting the optimal solution $\bm{x}^\ast$, all satellites have identical $a$, $e$, and $inc$ values but each satellite $i$ holds the following pair of $(\Omega_i, u_i)$: satellite 1 has $(\Omega_1,u_1)=(92.48\degree,105.12\degree)$, satellite 2 has $(\Omega_2,u_2)=(178.16\degree,311.04\degree)$, satellite 3 has $(\Omega_3,u_3)=(196.16\degree,203.04\degree)$, satellite 4 has $(\Omega_4,u_4)=(281.12\degree,53.28\degree)$, and satellite 5 has $(\Omega_5,u_5)=(6.80\degree,259.20\degree)$. Note that one can check that these satellites conform with the distribution rule shown in Eq.~\eqref{eq:raan_m} (by replacing $M_i$ with $u_i$ for circular orbits). Because $\bm{r}$ and $\bm{\pi}$ are time-invariant, we can easily approximate $Z_{\text{LP}}$ by computing Term~(1) of Eq.~\eqref{eq:zlpb}, resulting in $\hat{Z}_{\text{LP}}=410$. In fact, by directly solving the LP relaxation problem, we obtain $Z_{\text{LP}}=410$, which is identical to $\hat{Z}_{\text{LP}}$. The results are visualized in Fig.~\ref{fig:example3}.
\begin{figure}[htbp]
	\centering
	\begin{subfigure}[h]{0.45\linewidth}
		\centering
		\includegraphics[width=0.87\linewidth]{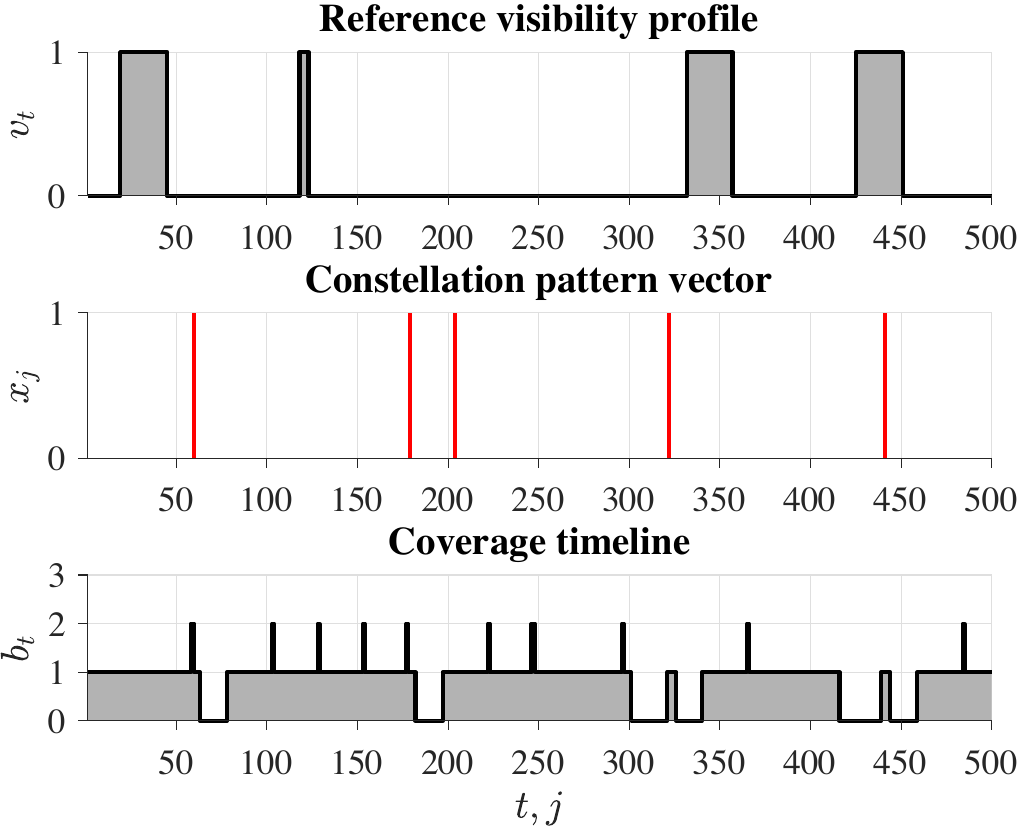}
		\caption{Illustration of $\bm{v}$, $\bm{x}$, and $\bm{b}$.}
	\end{subfigure}
	\begin{subfigure}[h]{0.45\linewidth}
		\centering
		\includegraphics[width=0.73\linewidth]{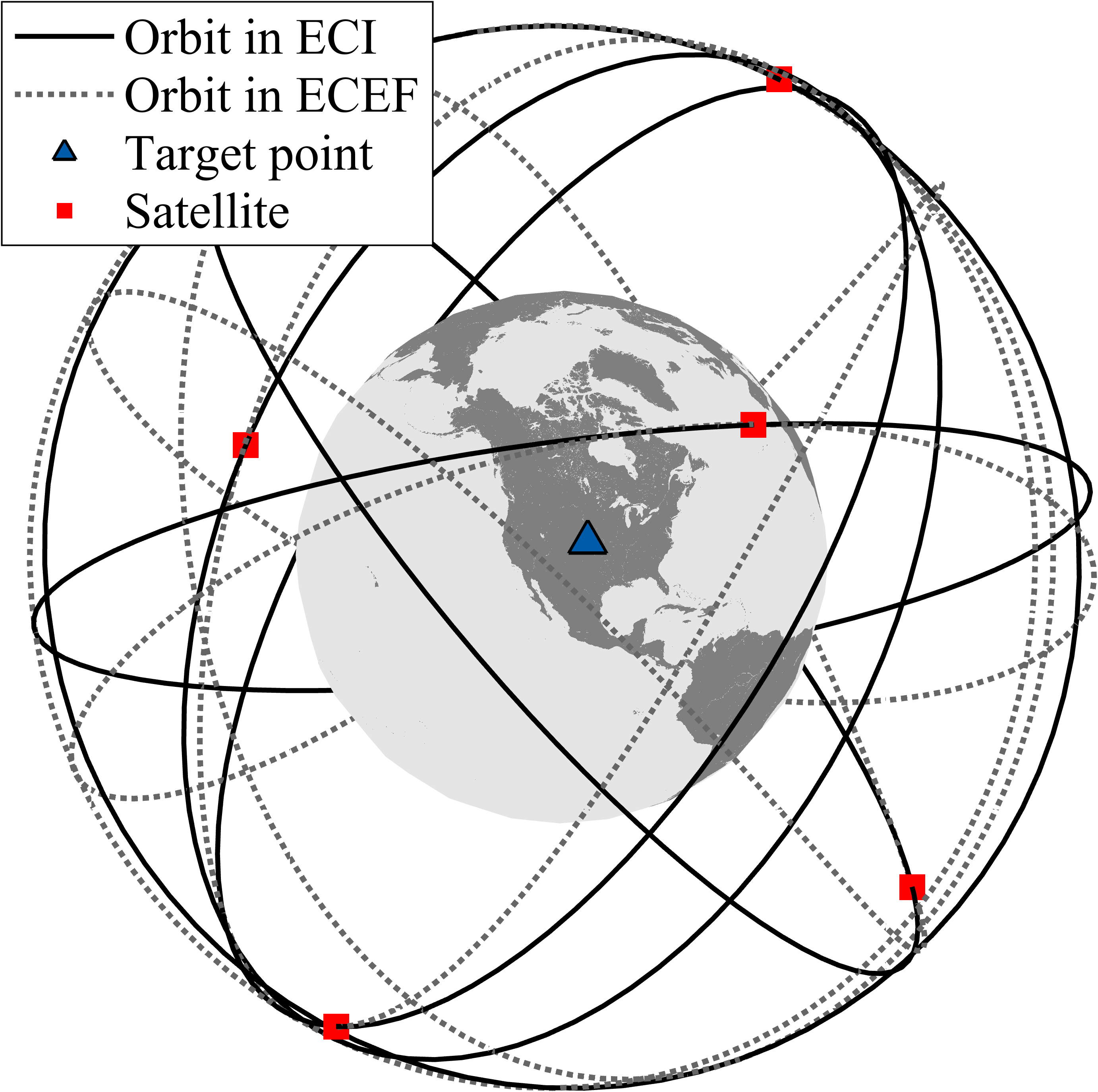}
		\caption{Corresponding constellation configuration in 3-D.}
	\end{subfigure}
	\caption{MCP solution for Example~\ref{eg:3}.}
	\label{fig:example3}
\end{figure}
\end{example}

\subsection{Constellation Transfer: Assignment Problem} \label{sec:ap}
The reconfiguration of a constellation incurs costs. In this subsection, we investigate the application of the Assignment Problem (AP) as a means to model the constellation transfer component of the reconfiguration process. Furthermore, we analyze the special mathematical properties of AP, which we subsequently utilize to devise a computationally efficient solution presented in Section~\ref{sec:solution}.

The transfer problem can be described over a bipartite graph $\mathcal{G}=(\mathcal{I} \cup \mathcal{J},\mathcal{E})$, in which the nodes of the set $\mathcal{I}$ describe the locations of the satellites while the nodes of the set $\mathcal{J}$ describe the locations of the orbital slots. Each edge $(i,j)\in\mathcal{E}$ is associated with the weight, or the cost $c_{ij}$ of transferring satellite $i$ to orbital slot $j$, commonly represented by $\Delta v$ required or the time-of-flight.

Within this framework, the transfer problem can be further divided into two main components: (i) the combinatorial optimization to find the minimum-cost assignment of satellites from one configuration to another and (ii) the orbital transfer trajectory design between a given (satellite, orbital slot) pair. The first component is concerned with the minimum-cost bipartite matching, which can be formulated as an assignment problem \cite{deweck2008optimal} as shown in Formulation~\ref{md:ap}. The second component deals with the construction of the cost matrix by evaluating the weights of edges in the bipartite graph setting. One can quantify each weight of the edges by solving an orbital boundary value problem. Because enumerating every edge can be time-consuming, several studies have proposed a rapid closed-form approximation of the true cost matrices \cite{avendano2009,davis2010constellation}. In this paper, we limit the scope of our work to high-thrust systems due to their benefit of timely reconfiguration, although low-thrust systems can also be considered as an alternative option.

\begin{formulation}[Assignment problem] \label{md:ap}
Let $\mathcal{I}=\{1,\dots,n\}$ denote the set of workers and $\mathcal{J}=\{1,\dots,m\}$ denote the set of projects. The cost of assigning worker $i$ to project $j$ is represented by $c_{ij}$. In the case of an unbalanced AP, the goal is to find the minimum-cost assignment of $n$ workers to $m$ projects such that all workers are assigned to projects, but not all projects are assigned with workers (i.e., when $n<m$). AP can be formulated as an integer linear program:
\begin{alignat*}{2}
	\text{(AP)} \quad \min \quad
	& \sum_{i\in \mathcal{I}} \sum_{j\in \mathcal{J}} c_{ij}\varphi_{ij} \\
	\text{s.t.} \quad
	& \sum_{j\in \mathcal{J}} \varphi_{ij}=1, \quad & \forall i\in \mathcal{I}\\
	& \sum_{i\in \mathcal{I}} \varphi_{ij}\le1, & \forall j\in \mathcal{J}\\
	& \varphi_{ij} \in \{0,1\}, & \forall i\in \mathcal{I},\forall j\in \mathcal{J}
\end{alignat*}
where the decision variable $\varphi_{ij}=1$ if worker $i$ is assigned to project $j$ ($\varphi_{ij}=0$ otherwise).
\end{formulation}

AP has garnered significant attention in the field of satellite constellation reconfiguration research as an optimization model for a constellation transfer problem. Introduced in a study by de Weck et al. \cite{deweck2008optimal}, a constellation transfer problem can be intuitively modeled as AP using the following analogy: the satellites as the workers and the orbital slots as the projects; the coefficient $c_{ij}$ as the cost (e.g., the fuel consumption) of transferring satellite $i$ to orbital slot $j$. The objective is to determine the minimum-cost assignment of $n$ satellites to $m$ orbital slots. In this paper, we adopt this transcription of the AP formulation in the modeling of a constellation transfer problem. For more information about the specific cost matrix generation used in this paper, refer to Ref.~\cite{vallado2013fundamentals}.

We briefly discuss the concepts of totally unimodular (TU) matrices and integral polyhedra, which are useful in leading to the discussion of the assignment problem. These concepts will come in handy later in this paper.

\begin{definition}[Total unimodularity] \label{def:tu}
    An integral matrix $\bm{A}$ is TU if every square sub-matrix of $\bm{A}$ has determinant equal to 0, 1, or -1.
\end{definition}

\begin{theorem}[Hoffman-Kruskal \cite{hoffman1956}] \label{thm:hk}
    Let $\bm{A}$ be an integral matrix. The polyhedron $\{\bm{x}:\bm{A}\bm{x}\le\bm{b},\bm{x}\ge \bm{0}\}$ is integral for all integral vector $\bm{b}$ if and only if $\bm{A}$ is TU.
\end{theorem}

One special feature of AP is that the problem satisfies Theorem~\ref{thm:hk}. That is, the constraint matrix of AP, also known as the incidence matrix of a bipartite graph, is totally unimodular and the right-hand vector is integral. As a result, the extreme points of the corresponding polytope are integral. We say that such a problem possesses the \textit{integrality property}. Consequently, the problem can be efficiently solved as a linear program (e.g., using the simplex or the interior-point methods) by relaxing the integrality constraints, also known as the linear programming relaxation, and still obtain integral optimal solutions. Other specialized algorithms such as the polynomial-time Hungarian algorithm \cite{kuhn1955} (also known as the Kuhn-Munkres algorithm; the asymptotic complexity is known to be $\mathcal{O}(m^3)$ for a square $m\times m$ matrix) or an auction algorithm (with pseudopolynomial complexity \cite{Bertsekas1981} and polynomial complexity using $\epsilon$-scaling \cite{Bertsekas1988}) are also available.

\section{Regional Constellation Reconfiguration Problem} \label{sec:integratedmodel}
\subsection{Problem Description}
Suppose a group of heterogeneous\footnote{The term heterogeneity embodies a general mission scenario of a federated system of satellites with different hardware specifications, orbital elements, and/or fuel states.} satellites is undertaking a reconfiguration process to form a new configuration to maximize the observation reward of a newly emerged set of spot targets, denoted as $\mathcal{P}$. Each target $p$ in $\mathcal{P}$ is associated with a time-dependent observation reward $\pi_{tp},\forall t \in\mathcal{T}$ where $\mathcal{T}$ is the set of time step indices. The reconfiguration process involves (i) the design of the maximum-reward destination configuration and (ii) the minimum-cost assignment of satellites between the initial and destination configurations. The goal of the problem is to identify a set of non-dominated solutions in the objective space spanned by these two competing objectives (i) and (ii).

We will refer to this problem as the Regional Constellation Reconfiguration Problem or RCRP in short. We use the term ``regional constellations'' to emphasize heterogeneity and asymmetry as distinctive design philosophies, in contrast to the homogeneity and symmetry of traditional global constellations.

\subsection{Mathematical Formulation} \label{sec:formulation}
In this subsection, we propose a mathematical formulation of the RCRP. In the proposed formulation, we consider a general case that accommodates multiple target points and multiple subconstellations \cite{lee2020satellite}. A subconstellation is a group of satellites in a given constellation system that share common orbital characteristics such as the semi-major axis, eccentricity, inclination, and argument of periapsis. Consequently, a constellation system can therefore consist of multiple subconstellations (an example would be a multi-layered constellation). 

Let $\mathcal{I}$ be the set of satellite indices (index $i$), $\mathcal{J}$ be the set of orbital slot indices (index $j$), $\mathcal{P}$ be the set of target point indices (index $p$), $\mathcal{S}$ be the set of subconstellation indices (index $s$), and    $\mathcal{T}$ be the set of time step indices (index $t$). By extending the concept of subconstellations to orbital slots, the set of all orbital slot indices $\mathcal{J}$ can be partitioned into $|\mathcal{S}|$ subsets such that $\mathcal{J}=\bigcup_{s\in\mathcal{S}}\mathcal{J}_s$ where $\mathcal{J}_s \subseteq \mathcal{J}$ denotes the set of orbital slot indices of subconstellation $s\in\mathcal{S}$, and $\mathcal{S}$ is the index set of $\mathcal{J}_s$. For the case with the RGT orbit assumption, we enforce the \textit{synchronous condition} to guarantee identical periods of repetition for all subconstellations: $T_s = T, \forall s \in \mathcal{S}$ where $T_s$ is the period of repetition for subconstellation $s$.

We also define the parameters $c_{ijs}$ as the cost of assigning satellite $i$ to orbital slot $j$ of subconstellation $s$ ($c_{ijs}\ge0$), $\pi_{tp}$ as the reward for covering target point $p$ at time step $t$ ($\pi_{tp}\ge0$), $r_{tp}$ as the coverage threshold for target point $p$ at time step $t$ ($r_{tp}\ge1$), and $b_{tp}$ as the number of satellite(s) in view from target point $p$ at time step $t$ ($b_{tp}\ge0)$. Furthermore, we define the parameter $V_{tjps}$ as follows:
\begin{equation*}
    V_{tjps} =\begin{cases}
      1, & \text{if a satellite in orbital slot $j$ of subconstellation $s$ covers target point $p$ at time step $t$} \\
      0, & \text{otherwise}\\
    \end{cases} 
\end{equation*}
The decision variables of interest include $\varphi_{ijs}$ and $y_{tp}$, defined as follows, respectively:
\begin{equation*}
     \varphi_{ijs} =\begin{cases}
1, &\text{if satellite $i$ is allocated to orbital slot $j$ of subconstellation $s$} \\
0, &\text{otherwise}
\end{cases} 
\end{equation*}
\begin{equation*}
    y_{tp}=\begin{cases}
1, &\text{if target point $p$ is covered at time step $t$ ($b_{tp}\ge r_{tp}$)} \\
0, &\text{otherwise}
\end{cases}
\end{equation*}
   
With the above notations, the mathematical formulation of RCRP is as follows
\begin{subequations}
\begin{alignat}{2}
\text{(RCRP)} \quad \min \quad & \left[ \begin{aligned}&\sum_{s\in \mathcal{S}}\sum_{j\in \mathcal{J}_s}\sum_{i\in \mathcal{I}}c_{ijs} \varphi_{ijs} \\
&-\sum_{p\in \mathcal{P}}\sum_{t\in \mathcal{T}} \pi_{tp} y_{tp}
\end{aligned} \right] \label{eq:j}\\
\text{s.t.} \quad & \sum_{s \in \mathcal{S}}\sum_{j\in \mathcal{J}_s}\varphi_{ijs} = 1, & \forall i\in \mathcal{I} \label{eq:c1}\\
& \sum_{i\in \mathcal{I}}\varphi_{ijs} \le 1, & \forall j\in \mathcal{J}_s,\forall s \in \mathcal{S}\label{eq:c2}\\
& \sum_{s\in \mathcal{S}}\sum_{j\in \mathcal{J}_s}\sum_{i\in \mathcal{I}} V_{tjps}\varphi_{ijs} \ge r_{tp}y_{tp}, & \forall t\in \mathcal{T}, \forall p\in \mathcal{P} \label{eq:c3}\\
& \varphi_{ijs} \in \{0,1\}, & \forall i\in \mathcal{I},\forall j\in \mathcal{J}_s,\forall s\in \mathcal{S} \label{eq:x} \\
& y_{tp} \in \{0,1\}, & \forall t\in \mathcal{T}, \forall p\in \mathcal{P} \label{eq:y}
\end{alignat}
\end{subequations}

RCRP is formulated as a bi-objective ILP. The first objective function in \eqref{eq:j} minimizes the total cost of a constellation reconfiguration process, while the second objective function in \eqref{eq:j} maximizes the total reward earned by covering a set of target points. Constraints~\eqref{eq:c1} and \eqref{eq:c2} are the AP-based constraints; Constraints~\eqref{eq:c1} ensure that every satellite is assigned to an orbital slot and Constraints~\eqref{eq:c2} restrict at most one satellite to be occupied per orbital slot. Constraints~\eqref{eq:c3} are the MCP-based constraints; these constraints ensure that the target point $p$ is covered at time step $t$ only if there exists at least $r_{tp}$ satellite(s) in view. Note that the cardinality constraint of MCP [i.e., Constraint~\eqref{eq:mcpb} of Formulation~\ref{md:mcp}] is omitted in this formulation because it is implied by the satellite indices set $\mathcal{I}=\{1,\dots,n\}$ and the AP constraints. Constraints~\eqref{eq:x} and \eqref{eq:y} define the domains of decision variables.

\begin{figure}[htbp]
	\centering
		\includegraphics[width=0.42\linewidth]{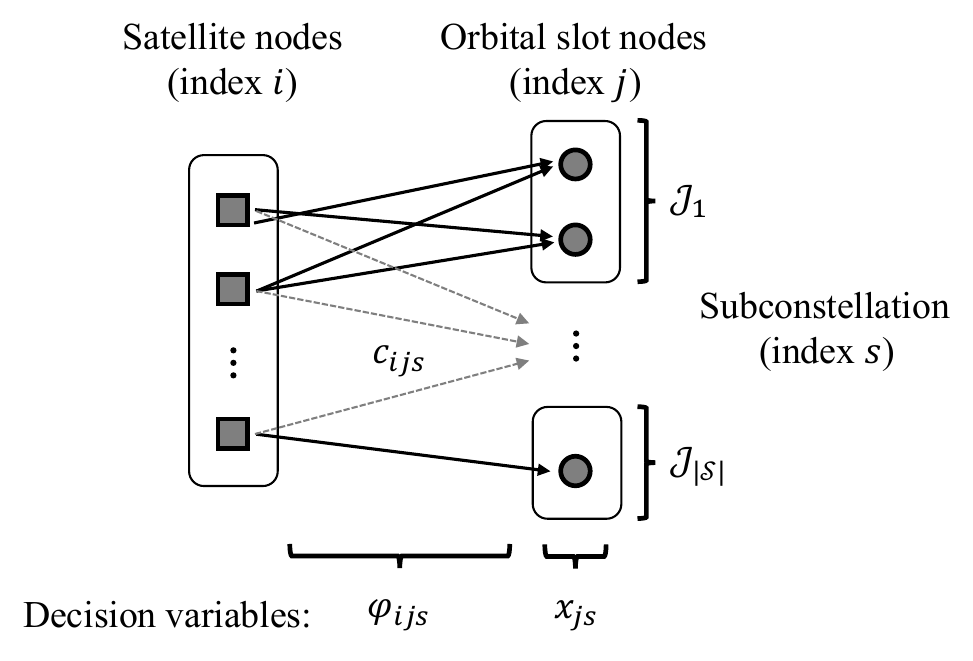}
		\caption{Decision variables of AP and MCP and their relationship.}
		\label{fig:variables}
\end{figure}

Notice the decision variables of RCRP---they are in the form of the AP decision variables; the reasoning behind this choice is explained. The decision variable $\varphi_{ijs}$ of AP indicates an assignment of satellite $i$ to orbital slot $j$ of subconstellation $s$, while the decision variable $x_{js}$ of MCP indicates whether a satellite occupies orbital slot $j$ of subconstellation $s$. Therefore, it follows naturally that $\varphi_{ijs}$ are the elemental decision variables because we can deduce $x_{js}$ from $\varphi_{ijs}$ (see Fig.~\ref{fig:variables}). The following relationship couples these two different sets of decision variables along with Constraints~\eqref{eq:c1} and \eqref{eq:c2}:
\begin{equation}
\label{eq:couple}
x_{js} = \sum_{i\in \mathcal{I}} \varphi_{ijs}, \quad \forall j\in \mathcal{J}, \forall s \in \mathcal{S}
\end{equation}
where both $\varphi_{ijs}$ and $x_{js}$ are binary variables. This coupled relationship in Eq.~\eqref{eq:couple} enables an integrated ILP formulation that simultaneously considers both the constellation transfer problem and the constellation design problem.

\subsection{Model Characteristics} \label{sec:charcteristics}

The RCRP formulation possesses the following characteristics. First, RCRP is NP-hard because of the embedded MCP structure (cf. Formulation~\ref{md:mcp}). This deduction follows from the NP-hardness of MCLP \cite{megiddo1983}, which has shown to be a particular case of MCP (see discussion in Section~\ref{sec:mcp}). Second, the AP structure [Constraints~\eqref{eq:c1},\eqref{eq:c2}, and \eqref{eq:x}] is preserved in RCRP with the decision variables being AP-based. In this perspective, the complicating constraints are Constraints~\eqref{eq:c3}.

The RCRP formulation combines the constellation transfer problem with the AP formulation and the constellation design problem with the MCP formulation. The former exhibits a special structure---the integrality property---that enables an efficient solution approach. The latter, however, is a combinatorial optimization problem, making the use of exact methods, such as the branch-and-bound algorithm, computationally expensive. In light of this observation, we develop a solution method in Section~\ref{sec:solution} that capitalizes on the characteristics of the RCRP formulation.

\section{Lagrangian Relaxation-Based Solution Method} \label{sec:solution}
This section develops a solution method for the RCRP, which is a bi-objective combinatorial optimization problem. To approach the bi-objective formulation, we use the $\varepsilon$-constraint method \cite{haimes1971} to transform RCRP into a single-objective optimization problem. This is then solved in series by varying the $\varepsilon$ value. The transformed single-objective problem can be solved using a commercial MILP solver; however, this approach can become computationally challenging even for moderately-sized instances. Motivated by this background and the need to rapidly characterize reconfiguration trade-offs for timely reconfiguration, we propose a computationally-efficient Lagrangian relaxation-based solution method that leverages the unique structure of the model.

\subsection{$\varepsilon$-constraint Reformulation}
The objective of RCRP is to identify a set of non-dominated solutions, as specified by its bi-objective formulation. To solve this problem, we reformulate RCRP as a single-objective optimization problem via the $\varepsilon$-constraint method by casting one of the two objective functions into a constraint with an upper (or lower) bound $\varepsilon$. Solving a series of single-objective $\varepsilon$-constrained problems to optimality given the sequence of $\varepsilon$ values yields the set of non-dominated solutions, also known as the Pareto front, of the original problem. Selecting an appropriate objective function for the constraint transformation is important because the choice made in this step affects the downstream algorithmic efforts.

Applying the $\varepsilon$-constraint method to RCRP, we transform the cost minimization objective function into a constraint that is bounded from above by $\varepsilon$ [Constraint~\eqref{eq:epsilon}]. In a physical sense, $\varepsilon$ represents the maximum allowable aggregated cost of reconfiguration, hence the name \textit{aggregated resource constraint} (ARC) for Constraint~\eqref{eq:epsilon}. The following is the single-objective model with the ARC:
\begin{align}
\text{(RCRP-ARC)} \quad Z(\varepsilon)= \min \quad & -\sum_{p\in \mathcal{P}}\sum_{t\in \mathcal{T}} \pi_{tp} y_{tp} \nonumber \\
\text{s.t.} \quad & \sum_{s\in \mathcal{S}}\sum_{j\in \mathcal{J}_s}\sum_{i\in \mathcal{I}}c_{ijs} \varphi_{ijs} \le \varepsilon \label{eq:epsilon} \\
& \text{Constraints~\eqref{eq:c1}--\eqref{eq:y}} \nonumber
\end{align}
where $Z(\varepsilon)$ denotes the optimal value of RCRP-ARC with a parameter of $\varepsilon$. The goal of RCRP-ARC is to maximize the total coverage reward while being subject to the aggregated resource constraint.

Algorithm~\ref{alg:overall} outlines the procedure based on the $\varepsilon$-constraint method, as described previously. It is important to note that the algorithm begins by solving an AP to determine the minimum-cost assignment of satellites for transitioning between configurations. The obtained optimum serves as a starting value, denoted as $\varepsilon_0$, for the $\varepsilon$-constraint method. If $\mathcal{I}\subset\mathcal{J}$, which indicates that satellites can maintain their orbits without executing any orbital maneuvers, the value of $\varepsilon_0$ is set to zero.

\subsection{Lagrangian Relaxation}
The Lagrangian relaxation is a decomposition-based optimization technique used to approach complex problems by dualizing complicating constraints, thereby exposing the remaining ``relatively easy'' structure for efficient solving (see Ref.~\cite{fisher2004} for a general overview of the topic). Specifically, in our case, the complicating constraints can be viewed as those of MCP, Constraints~\eqref{eq:c3}, primarily due to the intact AP structure (Section~\ref{sec:charcteristics}) present in the relaxed problem, along with Constraint~\eqref{eq:epsilon}. For a more in-depth discussion on selecting between competing relaxations, we invite readers to refer to Appendix~B.

To retrieve the \textit{Lagrangian problem} (LR) of RCRP-ARC, we dualize Constraints~\eqref{eq:c3}:\hypertarget{(LR)}{}
\begin{alignat}{2}
\text{(LR)} \quad Z_\text{D}(\varepsilon,\bm{\lambda}) = \min \quad & -\sum_{p\in \mathcal{P}}\sum_{t\in \mathcal{T}} \pi_{tp} y_{tp} + \sum_{p\in \mathcal{P}}\sum_{t\in \mathcal{T}} \lambda_{tp} \Bigg[r_{tp}y_{tp} - \sum_{s\in \mathcal{S}}\sum_{j\in \mathcal{J}_s}\sum_{i\in \mathcal{I}} V_{tjps}\varphi_{ijs}\Bigg] \nonumber \\
\text{s.t.} \quad & \sum_{s \in \mathcal{S}}\sum_{j\in \mathcal{J}_s}\varphi_{ijs} = 1, & \forall i\in \mathcal{I} \nonumber \\
& \sum_{i\in \mathcal{I}}\varphi_{ijs} \le 1, & \forall j\in \mathcal{J}_s,\forall s \in \mathcal{S} \nonumber \\
& \sum_{s\in \mathcal{S}}\sum_{j\in \mathcal{J}_s}\sum_{i\in \mathcal{I}}c_{ijs} \varphi_{ijs} \le \varepsilon \nonumber \\
& \varphi_{ijs} \in \{0,1\}, & \forall i\in \mathcal{I},\forall j\in \mathcal{J}_s,\forall s\in \mathcal{S} \nonumber \\
& y_{tp} \in \{0,1\}, & \forall t\in \mathcal{T},\forall p\in \mathcal{P} \nonumber
\end{alignat}
where $\bm{\lambda}=(\lambda_{tp} \in \mathbb{R}_{\geq0}:t \in \mathcal{T}, p \in \mathcal{P})$ is a vector of Lagrange multipliers associated with Constraints~\eqref{eq:c3}, and $Z_\text{D}(\varepsilon,\bm{\lambda})$ denotes the optimal value of LR.

\begin{remark}
\label{rm:lb}
For all non-negative $\bm{\lambda}$, we have $Z_\text{D}(\varepsilon,\bm{\lambda})\le Z(\varepsilon)$. It is easy to see this because for a given optimal solution $(\bm{\varphi}^\ast,\bm{y}^\ast)$ to RCRP-ARC, we observe that the following series of inequalities hold:
\begin{align*}
    Z(\varepsilon) & \ge -\sum_{p\in \mathcal{P}}\sum_{t\in \mathcal{T}} \pi_{tp} y_{tp}^\ast + \sum_{p\in \mathcal{P}}\sum_{t\in \mathcal{T}} \lambda_{tp} \Bigg[r_{tp}y_{tp}^\ast - \sum_{s\in \mathcal{S}}\sum_{j\in \mathcal{J}_s}\sum_{i\in \mathcal{I}} V_{tjps}\varphi_{ijs}^\ast\Bigg] \\
    & \ge Z_{\text{D}}(\varepsilon,\bm{\lambda})
\end{align*}
where the first inequality follows from adding the non-positive term to $Z(\varepsilon)$. The second inequality results from the fact that relaxing Constraint~\eqref{eq:c3} may expand the feasible region and potentially allow for the discovery of an improving solution that further reduces the objective function value.
\end{remark}

\subsection{Lagrangian Dual Problem and Subgradient Method}
Following from Remark~\ref{rm:lb}, we observe that the lower bound $Z_{\text{D}}(\varepsilon,\bm{\lambda})$ can be tightened up (i.e., maximized) by solving for the optimal $\bm{\lambda}^\ast$. Such a problem is called the \textit{Lagrangian dual problem} (D) and is formulated as follows:
\begin{equation*}
    \text{(D)} \quad Z_{\text{D}}(\varepsilon) = \max_{\bm{\lambda}} Z_{\text{D}}(\varepsilon,\bm{\lambda})
\end{equation*}

The Lagrangian dual problem is a non-differentiable optimization problem because $Z_{\text{D}}(\varepsilon,\bm{\lambda})$ is a piecewise linear, concave function of $\bm{\lambda}$. To solve this problem, we employ the subgradient method \cite{held1971}, which has been demonstrated as an effective method for non-differentiable optimization problems \cite{held1974}. The subgradient method is an iterative algorithm in the spirit of the gradient ascent method for determining the maximum solution of a continuously differentiable function. Algorithm~\ref{alg:subgradient} provides the pseudocode for the subgradient optimization.

The subgradient method starts by initializing the Lagrange multipliers, $\bm{\lambda}^0$. At iteration $k$, the Lagrangian problem $\text{LR}^k$ is solved using the provided parameters, $\varepsilon$, $\bm{\lambda}^k$, $\bm{c}$, $\bm{\pi}$, $\bm{r}$, and $\bm{v}$. With the optimal solution $(\bm{\varphi}^k,\bm{y}^k)$ to $\text{LR}^k$, we apply a local search-based heuristic method to obtain an estimate $\hat{Z}(\varepsilon)$ of the optimal value $Z(\varepsilon)$. To ensure the best estimate of the optimal value is maintained and to aid in the convergence of the subgradient method, the best $\hat{Z}(\varepsilon)$ up to iteration $k$ is stored in memory as the incumbent optimum.

Next, the subgradient $\bm{g}^k$ of $Z_{\text{D}}(\varepsilon,\bm{\lambda}^k)$ at $\bm{\lambda}^k$ is computed:
\begin{equation*}
    \label{eq:subgradient}
    \bm{g}^k=\Bigg(r_{tp}y_{tp}^k - \sum_{s\in \mathcal{S}}\sum_{j\in \mathcal{J}_s}\sum_{i\in \mathcal{I}} V_{tjps}\varphi_{ijs}^k:t\in\mathcal{T},p\in\mathcal{P}\Bigg)
\end{equation*}
If the subgradient of $Z_{\text{D}}(\varepsilon,\bm{\lambda}^k)$ at $\bm{\lambda}^k$ is $\bm{0}$, then $\bm{\lambda}^k$ is an optimal Lagrange multiplier vector, and the algorithm terminates. The algorithm may also terminate with suboptimal Lagrangian multipliers if any of the following termination criteria are triggered: the maximum iteration count, the gap tolerance between $\hat{Z}(\varepsilon)$ and $Z_{\text{D}}(\varepsilon,\bm{\lambda}^k)$, and the step size tolerance. With $\bm{\lambda}^0=\bm{0}$, the Lagrangian relaxation bound starts with $Z_{\text{D}}(\varepsilon,\bm{\lambda}^0)=-\sum_{p\in \mathcal{P}}\sum_{t\in \mathcal{T}} \pi_{tp}$ and improves as the subgradient method progresses. In the case of premature termination due to reaching the maximum iteration limit, the obtained $\bm{\lambda}$ may be suboptimal, resulting in $Z_{\text{D}}(\varepsilon,\bm{\lambda})< Z_{\text{D}}(\varepsilon)$. With the knowledge of the optimal dual variables of the LP relaxation problem, one can use them for $\bm{\lambda}^0$. However, this approach is not ideal as it requires running the LP relaxation problem beforehand, which can be computationally expensive for large instances.

Unless the termination flag is triggered, the algorithm reiterates the procedure with the new set of Lagrange multipliers. The rule for updating the Lagrange multipliers is as follows:
\begin{equation*}
\bm{\lambda}^{k+1} \coloneqq \max(\bm{0},\bm{\lambda}^k+\theta_k\bm{g}^k)
\end{equation*}
where $\max(\cdot)$ is the element-wise maximum to guarantee the non-negativity of $\lambda_{tp}$. Unless the dualized constraints are equality constraints, which can have associated multipliers unrestricted in sign, the multipliers need to be non-negative to penalize the violated constraints correctly \cite{bertsimas-LPbook}.

The step size $\theta_k$ commonly used in practice is:
\begin{equation}
    \label{eq:stepsize}
    \theta_k\coloneqq\frac{\hat{Z}(\varepsilon)-Z_{\text{D}}(\varepsilon,\bm{\lambda}^k)}{\lVert \bm{g}^k \rVert^2}\alpha_k
\end{equation}
where $\lVert \cdot \rVert$ is the Euclidean norm, and $\alpha_k$ is a scalar satisfying $0<\alpha_k\le 2$. The proof of convergence of the above step size formula is referred to Ref.~\cite{held1974}. As recommended by Fisher \cite{fisher2004}, the starting value of $\alpha_k$ is set to $\alpha_0=2$ and is halved if $Z_{\text{D}}(\varepsilon,\bm{\lambda})$ fails to increase in a number of iterations. While there are different types of $\theta_k$ proposed in literature, the step size formula in Eq.~\eqref{eq:stepsize} has performed particularly well in our problem settings.

The subgradient method suffers from several drawbacks, such as the zigzagging phenomenon and slow convergence to the optimal multipliers $\bm{\lambda}^\ast$. Several studies have proposed variants of the subgradient method, such as the surrogate Lagrangian relaxation method and the bundle method, to alleviate these issues. Interested readers can refer to Ref.~\cite{guignard2013} for additional materials on methods for non-differentiable problems.

There exist two computational bottlenecks in this algorithm. One at computing the lower bound $Z_{\text{D}}(\varepsilon,\bm{\lambda}^k)$ and another at computing the upper bound $\hat{Z}(\varepsilon)$, both of which occur at every iteration. To make each iteration of the subgradient method more efficient, we provide efficient ways to compute $Z_{\text{D}}(\varepsilon,\bm{\lambda}^k)$ and $\hat{Z}(\varepsilon)$.

\subsection{Lower Bound: Lagrangian Problem Decomposition} \label{sec:lb}
At each iteration of the subgradient optimization, $Z_{\text{D}}(\varepsilon,\bm{\lambda}^k)$ is computed. By relaxing the complicating constraints [Constraints~\eqref{eq:c3}], which are also the linking constraints, we observe that Problem~\LR can be decomposed into two subproblems based on the variable types, $\bm{\varphi}$ and $\bm{y}$: 
\LRo, an assignment problem with a side constraint, and \LRt, an unconstrained binary integer linear program. \hypertarget{(LR1)}{}\begin{alignat*}{2}
(\text{LR}1) \quad Z_\text{D1}(\varepsilon,\bm{\lambda}) = \min \quad & - \sum_{p\in \mathcal{P}}\sum_{t\in \mathcal{T}} \lambda_{tp} \Bigg[\sum_{s\in \mathcal{S}}\sum_{j\in \mathcal{J}_s}\sum_{i\in \mathcal{I}} V_{tjps}\varphi_{ijs}\Bigg] \\
\text{s.t.} \quad & \sum_{s \in \mathcal{S}}\sum_{j\in \mathcal{J}_s}\varphi_{ijs} = 1, & \forall i\in \mathcal{I}  \\
& \sum_{i\in \mathcal{I}}\varphi_{ijs} \le 1, & \forall j\in \mathcal{J}_s,\forall s \in \mathcal{S} \\
& \sum_{s\in \mathcal{S}}\sum_{j\in \mathcal{J}_s}\sum_{i\in \mathcal{I}}c_{ijs} \varphi_{ijs} \le \varepsilon \\
& \varphi_{ijs} \in \{0,1\}, & \forall i\in \mathcal{I},\forall j\in \mathcal{J}_s,\forall s\in \mathcal{S}
\end{alignat*}
\hypertarget{(LR2)}{}\begin{alignat*}{2}
(\text{LR}2) \quad Z_\text{D2}(\bm{\lambda}) = \min \quad & \sum_{p\in \mathcal{P}}\sum_{t\in \mathcal{T}} \big(\lambda_{tp} r_{tp}-\pi_{tp} \big)y_{tp} \\
\text{s.t.} \quad & y_{tp} \in \{0,1\}, & \forall t\in \mathcal{T},\forall p\in \mathcal{P}
\end{alignat*}

For convenience, we denote $\Phi \coloneqq\{\bm{\varphi}\in\{0,1\}^{|\mathcal{I}||\mathcal{J}|}:\text{Constraints~\eqref{eq:c1}, \eqref{eq:c2}, \eqref{eq:epsilon}}\}$ the set of assignments $\bm{\varphi}$ satisfying the AP and ARC constraints.

Problem~\LR is characterized as follows:
\begin{enumerate}[label=(\roman*)]
    \item Let $\varepsilon_{\text{r}}$ denote the critical value of $\varepsilon$ at which Constraint~\eqref{eq:epsilon} becomes redundant in $\text{Co}(\Phi)$, the convex hull of the feasible region of Problem~\LRo. From inspection, it is evident that no total cost of assignments exceeds the value of
    \begin{equation*}
        \varepsilon_{\text{r}} \coloneqq \sum_{i \in \mathcal{I}} \max_{j\in\mathcal{J}_s,s\in\mathcal{S}} c_{ijs}
    \end{equation*}
    If $\varepsilon\ge\varepsilon_{\text{r}}$, then Constraint~\eqref{eq:epsilon} is redundant and can be removed without altering the convex hull. This reveals the intact AP structure [Constraints~\eqref{eq:c1},\eqref{eq:c2}, and \eqref{eq:x}], which possesses the integrality property (cf. Formulation~\ref{md:ap}). In this case, Problem~\LRo can be efficiently solved using a specialized AP algorithm or through LP. Otherwise, Problem~\LRo does not possess the integrality property and cannot be solved by means of LP. Based on our experience, the problem can be considered a relatively easy ILP and thus can be expected to be solved in a reasonable amount of time in many instances, albeit there may be instances where guaranteeing the solution optimality takes longer.
    \item Problem~\LRt can be solved trivially. For fixed multipliers, the optimal solution is generated by inspecting the coefficients of the objective function:
\begin{equation*}
    y_{tp}^\ast = \begin{cases}
      1, & \text{if } (\lambda_{tp} r_{tp}-\pi_{tp})<0 \\
      0, & \text{otherwise}
    \end{cases}
\end{equation*}
    \item The optimal value of Problem~\LR is the sum of the optimal values of its subproblems: $Z_{\text{D}}(\varepsilon,\bm{\lambda})=Z_{\text{D1}}(\varepsilon,\bm{\lambda})+Z_{\text{D2}}(\bm{\lambda})$.
\end{enumerate}

Using the argument of LP duality \cite{geoffrion1974}, one can deduce that the Lagrangian relaxation bound is at least as good as the LP relaxation bound: $Z_{\text{LP}}(\varepsilon) \le Z_{\text{D}}(\varepsilon)$. Based on observations made about Problem~\LR, we provide a discussion on the tightness of these two bounds and a sufficient condition at which they are equal.

The relative strength of $Z_{\text{D}}(\varepsilon)$ with respect to $Z_{\text{LP}}(\varepsilon)$ depends on the redundancy of Constraint~\eqref{eq:epsilon} in Problem LR1:
\begin{enumerate}[label=(\roman*)]
    \item If Constraint~\eqref{eq:epsilon} is not redundant, the Lagrangian relaxation may provide a tighter bound than the LP relaxation bound: $Z_{\text{LP}}(\varepsilon) \le Z_\text{D}(\varepsilon)$. The value of $\varepsilon$ dictates the tightness of these bounds.
    \item If Constraint~\eqref{eq:epsilon} is redundant, Problem~\LR possesses the integrality property, which is a sufficient condition for $Z_{\text{LP}}(\varepsilon) = Z_\text{D}(\varepsilon)$ \cite{geoffrion1974}.
\end{enumerate}

Reiterating, Constraint~\eqref{eq:epsilon} being redundant is only a sufficient condition for the equal bounds, and therefore it does not necessarily imply $Z_{\text{LP}}(\varepsilon)<Z_{\text{D}}(\varepsilon)$ for all $\varepsilon<\varepsilon_{\text{r}}$. In some instances, we observe that there exists $\varepsilon_{\text{c}}<\varepsilon_{\text{r}}$ such that $Z_{\text{LP}}(\varepsilon)<Z_{\text{D}}(\varepsilon)$ for $\varepsilon<\varepsilon_{\text{c}}$ (and $Z_{\text{LP}}(\varepsilon)=Z_{\text{D}}(\varepsilon)$ for $\varepsilon\ge\varepsilon_{\text{c}}$). Note that if $Z_{\text{LP}}=Z_{\text{D}}$, then the optimal Lagrange multipliers are equal to the dual variables associated with Constraints~\eqref{eq:c3} in the LP relaxation of RCRP-ARC.

Revisit the earlier discussion on the LP relaxation bound for MCP discussed in Section~\ref{sec:mcp}. There is a special case, that is, $\xi_t=\xi, \forall t \in \mathcal{T}$ and $\bm{V}$ being a circulant matrix, for which the value of $\hat{Z}_{\text{LP}}$ is independent of optimal LP solution $\bm{x}^\ast$. This analysis is readily extensible to the RCRP-ARC formulation. Consequently, for $\varepsilon\ge\varepsilon_{\text{r}}$, $\xi_{tp}=\xi_p,\forall t \in \mathcal{T}$, and all relevant $\bm{V}$ matrices being circulant, the approximated LP relaxation bound for RCRP-ARC is simply [cf. Eq.~\eqref{eq:zlpb}]:
\begin{equation*}
    \hat{Z}_{\text{LP}}\coloneqq \max\bigg(-|\mathcal{I}|\sum_{p\in\mathcal{P}}\sum_{t\in\mathcal{T}}\xi_{tp}v_{tp},-\sum_{p\in\mathcal{P}}\sum_{t\in\mathcal{T}} \pi_{tp}\bigg)
\end{equation*}
where $\xi_{tp}\coloneqq \pi_{tp}/r_{tp}$.

The value of $\hat{Z}_{\text{LP}}$ provides a conservative lower bound for both $Z_{\text{LP}}$ and $Z_{\text{D}}$ that requires no knowledge of the optimal LP solution $\bm{\varphi}^\ast$. As will be discussed later in the computational experiments in Section~\ref{sec:experiments}, we can see that $\hat{Z}_{\text{LP}}$ outperforms $Z_{\text{D}}$ for suboptimal Lagrange multipliers. However, its usefulness is limited to instances triggering the first argument of the max function.

\subsection{Upper Bound: Lagrangian Heuristic via Combinatorial Neighborhood Local Search} \label{sec:ub}
To evaluate the duality gap over the convergence (e.g., to be used for the definition of the step size), we need a method to find a feasible solution to the primal problem (RCRP-ARC) at each iteration $k$ of the subgradient method. While it is theoretically possible for a solution to $\text{LR}^k$ to be discovered as feasible to the primal problem in the course of solving the Lagrangian dual problem, this happens rarely \cite{fisher2004}. Instead, we attempt to convert the solution to $\text{LR}^k$ into a feasible solution to the primal problem. The reason for solution infeasibility is due to the inconsistency between weakly coupled $\bm{\varphi}^k$ and $\bm{y}^k$ solutions; the coupling comes from the fact that Problems $\text{LR1}^k$ and $\text{LR2}^k$ are linked only through $\bm{\lambda}^k$ in their objective functions. Intuitively speaking, the coverage state of the obtained constellation configuration, which is computed from the assignments using Eq.~\eqref{eq:ytp}, does not match the obtained coverage state. The key idea in finding a feasible solution is to drive the solution obtained at each iteration of the subgradient method to feasibility by employing a heuristic approach. Such an approach that exploits a solution produced in the course of solving the Lagrangian dual problem is called the Lagrangian heuristic in literature.

A simple and straightforward way is to accept the subgradient assignment $\bm{\varphi}^k$ as the valid primal solution and obtain the conforming coverage state $\tilde{\bm{y}}^k(\bm{\varphi}^k)$ (note that we distinguish it from $\bm{y}^k$, without a tilde), which is simply the coverage state of the constellation configuration obtained from the set of assignments $\bm{\varphi}^k$:
\begin{equation}
\label{eq:ytp}
    \tilde{y}_{tp}^k(\varphi_{ijs}^k)=\begin{cases}
    1, & \text{if } \displaystyle \sum_{s\in \mathcal{S}}\sum_{j\in \mathcal{J}_s}\sum_{i\in \mathcal{I}} V_{tjps}\varphi_{ijs}^k \ge r_{tp} \\
    0, & \text{otherwise}
    \end{cases}
\end{equation}
This approach always yields feasible solutions to the primal problem because it circumvents the inconsistency between $\bm{\varphi}^k$ and $\bm{y}^k$. The other way around, obtaining $\tilde{\bm{\varphi}}^k(\bm{y}^k)$ from a given $\bm{y}^k$, is infeasible in most cases or requires a combinatorial optimization approach in otherwise rarely feasible cases.

A more sophisticated approach is to employ a local search in the neighborhood of $\bm{\varphi}^k$ to improve the quality of the initial primal solution at the cost of additional computation time. In this paper, we propose the following definition of the \textit{1-exchange neighborhood}:
\begin{equation}
\label{eq:neighborhood}
    \mathcal{N}(\bm{\varphi})\coloneqq\{\bm{\varphi}^\prime\in\Phi:\text{$\bm{\varphi}^\prime$ obtained from $\bm{\varphi}$ by \textit{exchanging at most one} (satellite, orbital slot) pair}\}
\end{equation}

\begin{figure}[htbp]
	\centering
		\includegraphics[width=0.37\linewidth]{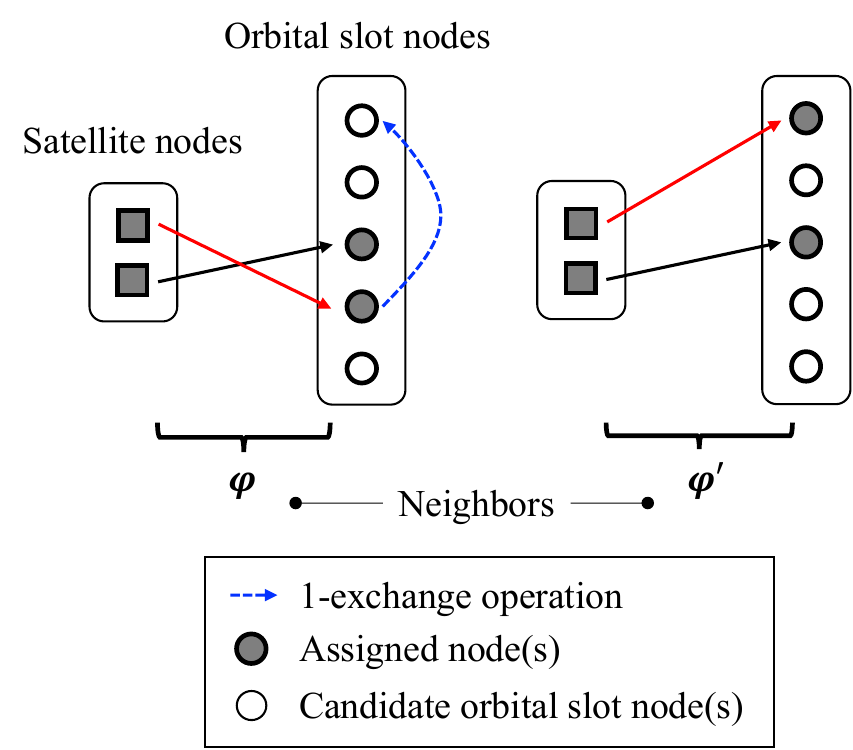}
		\caption{Illustration of the 1-exchange operation.}
		\label{fig:exchange}
\end{figure}

Figure~\ref{fig:exchange} illustrates the 1-exchange operation. Two sets of assignments, $\bm{\varphi}$ and $\bm{\varphi}'$, are \textit{neighbors} because at most a single (satellite, orbital slot) assignment pair is different, as indicated by the dashed line in the figure. Orbital slots that are occupied by other satellites are not part of the neighborhood because this would violate the definition of the feasible set $\Phi$. Additionally, not all unfilled orbital slots are valid candidates because all candidate orbital slots must conform with Constraint~\eqref{eq:epsilon}.

The 1-exchange neighborhood local search algorithm is presented in Algorithm~\ref{alg:heuristic}. Given $\bm{\varphi}$, the 1-exchange neighborhood $\mathcal{N}(\bm{\varphi})$ is defined as in Eq.~\eqref{eq:neighborhood}. The best candidate solution $\bm{\varphi}^\ast$ is found by exhaustively evaluating all candidate solutions in $\mathcal{N}(\bm{\varphi})$. Here, the figure of merit used is $\hat{Z}(\varepsilon)$. If the best candidate solution $\bm{\varphi}^\ast$ outperforms the incumbent solution, it is accepted as the new incumbent solution, and $\hat{Z}(\varepsilon)$ and the neighborhood are updated accordingly. The search process then reiterates with the new neighborhood. If not, the local search halts and the local optimum is obtained. The algorithm returns the local optimal solution $\bm{\varphi}^\ast$, its conforming $\tilde{\bm{y}}^\ast(\bm{\varphi}^\ast)$ [Eq.~\eqref{eq:ytp}], and the corresponding upper bound value $\hat{Z}(\varepsilon)$.

The trade-off between solution quality and computational effort exists in Line~3 of Algorithm~\ref{alg:heuristic} at deciding the local optimum of the current neighborhood. A brute-force evaluation of all candidate solutions can be costly for large instances. As such, it is practical to search only within a subset of the neighborhood $\mathcal{N}^\prime \subseteq \mathcal{N}$, which can be either randomly generated or defined by a pre-determined rule. Moreover, a first-come-first-served scheme can be applied to select the first local solution that improves the incumbent solution, effectively reducing the dimension of the search space. It is important to note that the radius of the neighborhood can substantially influence this trade-off. Employing a more general $\kappa$-exchange neighborhood local search ($\kappa>1$) can improve the quality of the heuristic solution for large instances, but at the cost of increased computational effort.

\subsection{Extension: RCRP with Individual Resource Constraints}\label{sec:extensions}
Various extensions can be made to the proposed formulation and method. One practically important example is discussed here. Constraint~\eqref{eq:epsilon} in the RCRP-ARC formulation limits the aggregated resource consumed by all system satellites. Along a similar line, but instead, we can formulate a variant of RCRP to enforce the maximum individual resource consumption on each satellite. This modeling approach has significant practical implications, as not all satellites have identical fuel states prior to a reconfiguration. For example, one can envision a realistic scenario of bringing a group of satellites with different fuel states together to form a federation for a new Earth observation mission \cite{golkar2015}.

We formulate the RCRP with \textit{individual resource constraints} (RCRP-IRC) as follows:
\begin{alignat}{2}
\text{(RCRP-IRC)} \ \ \min \quad & -\sum_{p\in \mathcal{P}}\sum_{t\in \mathcal{T}} \pi_{tp} y_{tp} \nonumber\\
\text{s.t.} \quad & \sum_{s\in \mathcal{S}}\sum_{j\in\mathcal{J}_s}c_{ijs}\varphi_{ijs} \le \varepsilon_i, \quad \forall i \in \mathcal{I}^\prime \label{eq:aggdeltav} \\
& \text{Constraints~\eqref{eq:c1}--\eqref{eq:y}} \nonumber
\end{alignat}
where $\mathcal{I}^\prime \subseteq \mathcal{I}$ represents the subset of satellites with IRC; Constraints~\eqref{eq:aggdeltav} define the reachable domain of orbital slots by satellite $i$ given the maximum allowance $\varepsilon_i$.

The proposed solution procedure based on the Lagrangian relaxation method for RCRP-ARC is readily applicable to RCRP-IRC. This is because the Lagrangian problem for RCRP-IRC is also separable into two subproblems based on the variable type, and the primal heuristic can be readily applied with the modified definition of $\Phi$. It should be noted that the analyses presented in Section~\ref{sec:lb} can be extended to RCRP-IRC with the appropriate values of $\varepsilon_{\text{r},i}$.

\section{Computational Experiments} \label{sec:experiments}
We conduct computational experiments to evaluate the performance of the proposed Lagrangian relaxation-based solution method. In particular, we focus on analyzing the solution quality and the computational efficiency of the Lagrangian heuristic in comparison to the results obtained by a mixed-integer programming (MIP) solver. We first perform the design of experiments in Section~\ref{sec:setup} and then compare the results obtained by the Lagrangian heuristic and a commercial software package in Section~\ref{sec:results}. The primary computational experiments are performed using RCRP-ARC for RGT orbits. In Section~\ref{sec:illustrative_example}, we provide an illustrative example to demonstrate the versatility of the proposed framework by extending it to a more general case of non-RGT orbits and RCRP-IRC.

\subsection{Test Instances} \label{sec:setup}

We generate test instances by varying the cardinalities of the sets $\mathcal{I}$, $\mathcal{J}$, $\mathcal{T}$, and $\mathcal{P}$. Each test instance selects one value from each of the following sets: $|\mathcal{I}| \in \{10,20\}$, $|\mathcal{J}|\in\{500,1000,2000\}$, and $|\mathcal{P}|\in\{10,20,30\}$. We assume, without loss of generality, that $|\mathcal{T}|=|\mathcal{J}|$.

Table~\ref{tab:size} presents the sizes of 18 randomly generated test instances for the RCRP problem. The instance pool ranges from the smallest, which contains up to \num{8.9e26} potentially feasible reconfiguration processes, to the largest, which contains up to \num{9.5e65} potentially feasible reconfiguration processes. Note that in Table~\ref{tab:size}, the constraint counts exclude the domain definitions of decision variables.

\begin{table}[htbp]
    \fontsize{9}{10}\selectfont
	\renewcommand{\arraystretch}{1.2}
	\caption{Test instances for RCRP-ARC and their sizes.}
	\centering
	\begin{tabular}{@{}r r r r r r r r@{}}
	\hline
	\hline
	\multirow{2}{*}{Instance} & \multirow{2}{*}{$|I|$} & \multirow{2}{*}{$|J|,|T|$} & \multirow{2}{*}{$|P|$} & \multicolumn{4}{c}{Size of RCRP-ARC} \\
	\cmidrule(lr){5-8}
	& & & & $\bm{\varphi}$ variables & $\bm{y}$ variables & Total variables & Total constraints \\
	\hline
1 & 10 & 500 & 10 & 5,000 & 5,000 & 10,000 & 5,511 \\
2 & 20 & 500 & 10 & 10,000 & 5,000 & 15,000 & 5,521 \\
3 & 10 & 500 & 20 & 5,000 & 10,000 & 15,000 & 10,511 \\
4 & 20 & 500 & 20 & 10,000 & 10,000 & 20,000 & 10,521 \\
5 & 10 & 1,000 & 10 & 10,000 & 10,000 & 20,000 & 11,011 \\
\arrayrulecolor{gray!50}\hline
6 & 20 & 1,000 & 10 & 20,000 & 10,000 & 30,000 & 11,021 \\
7 & 10 & 500 & 30 & 5,000 & 15,000 & 20,000 & 15,511 \\
8 & 20 & 500 & 30 & 10,000 & 15,000 & 25,000 & 15,521 \\
9 & 10 & 1,000 & 20 & 10,000 & 20,000 & 30,000 & 21,011 \\
10 & 20 & 1,000 & 20 & 20,000 & 20,000 & 40,000 & 21,021 \\
\hline
11 & 10 & 2,000 & 10 & 20,000 & 20,000 & 40,000 & 22,011 \\
12 & 20 & 2,000 & 10 & 40,000 & 20,000 & 60,000 & 22,021 \\
13 & 10 & 1,000 & 30 & 10,000 & 30,000 & 40,000 & 31,011 \\
14 & 20 & 1,000 & 30 & 20,000 & 30,000 & 50,000 & 31,021 \\
15 & 10 & 2,000 & 20 & 20,000 & 40,000 & 60,000 & 42,011 \\
\hline
16 & 20 & 2,000 & 20 & 40,000 & 40,000 & 80,000 & 42,021 \\
17 & 10 & 2,000 & 30 & 20,000 & 60,000 & 80,000 & 62,011 \\
18 & 20 & 2,000 & 30 & 40,000 & 60,000 & 100,000 & 62,021 \\
	\arrayrulecolor{black}\hline
	\hline
	\end{tabular}
	\label{tab:size}
\end{table}

For each RCRP test instance, we generate 10 RCRP-ARC sub-instances with varying values of $\varepsilon$. Without loss of generality, we set $\varepsilon_{\max}=\max c_{ijs}$ and create a sequence of 10 steps in the interval $[0,\varepsilon_{\max}]$. In total, we evaluate 18 instances for RCRP, which is equivalent to assessing 180 instances of RCRP-ARC.

Our goal is to capture a wide spectrum of orbital characteristics of orbital slots and $\vartheta_{\min}$. To do so, we randomly generate 18 parameter sets from the parameter space $\{N_{\text{P}} \in \mathbb{Z}_{\geq0}:30 \le N_{\text{P}} \le 45 \} \times \{inc \in \mathbb{R}_{\geq0}:0\degree \le inc \le 120\degree \} \times \{ \vartheta_{\min} \in \mathbb{R}_{\geq0}:5\degree \le \vartheta_{\min} \le 20 \degree\}$. Therefore, each RCRP instance has a unique set of $N_{\text{P}}$, $inc$, and $\vartheta_{\min}$ values. For all other parameters, the following values or generation rules are used. For the common orbital characteristics of $\mathcal{J}$, we let $N_{\text{D}}=3$ and $e=0$ (circular orbits). Each orbital slot $j\in\mathcal{J}$ has a pair of $(\Omega_j,u_j)$ that satisfies the distribution rule of the RGT common ground track constellation \cite{lee2020satellite}: $N_{\text{P}} (\Omega_j-\Omega_0)+N_{\text{D}} (u_j-u_0) = 0 \ \text{mod} \ 2\pi$ where we let $\Omega_0 = 0\degree$, $u_0=0\degree$. We assume $r_{tp}=1$ for all $t$ and $p$. The cost matrix is produced using the combined plane change and the Hohmann transfer maneuvers, as well as phasing maneuvers \cite{vallado2013fundamentals}. Additionally, we randomly generate the initial positions of $|\mathcal{I}|$ satellites in $\mathcal{J}$ from a discrete uniform distribution between 0 and $|\mathcal{J}|-1$. The target points are randomly distributed globally, bounded latitudinally by the inclination of a given test instance; the longitude and latitude each take a value from a uniform distribution. In summary, we can characterize the orbital slots as circular low Earth orbits covering a wide spectrum of inclination from prograde to retrograde (specifically, the altitude ranges between \SI{478.86}{km} and \SI{2729.95}{km}) and the minimum elevation angle threshold ranging from \SI{5}{\degree} to \SI{20}{\degree}.

We compare the results of the Lagrangian heuristic and a commercial MIP solver, Gurobi optimizer 9.1.1. Gurobi, the state-of-the-art solver for MILP problems, is chosen as the benchmark because there is no known specialized solver for RCRP. Note that existing algorithms for constellation reconfiguration are not suitable for this purpose as they do not account for the added layer of MCP constraints present in the RCRP formulation. Thus, to establish a benchmark for the performance of the Lagrangian heuristic, it is judicious to compare it against a general-purpose but widely-used MIP solver such as the Gurobi optimizer. The Gurobi optimizer utilizes an array of MIP techniques, including but not limited to, presolve, branch-and-bound, cutting plane, heuristics, and parallelism, at various phases of optimization to enhance the optimization efficiency with regard to both computation runtime and solution quality.

All computational experiments are coded in MATLAB and executed on a platform with an Intel Core i7-9700 3.00 GHz CPU processor (8 cores and 8 threads) and 32 GB memory. In all cases, we let the Gurobi optimizer utilize all 8 cores. The default Gurobi parameters are used except for the duality gap tolerance of \SI{0.5}{\%} (for both the baseline Gurobi case and \LRo) and the runtime limit of \SI{3600}{s}. If the Gurobi optimizer has not converged within the runtime limit, it returns the best incumbent primal solution found thus far. For the Lagrangian heuristic, we limit the 1-exchange neighborhood local search with the size $|\mathcal{N}^\prime|\le 10|\mathcal{I}|$ for the primal heuristic and the Gurobi optimizer for solving Problem~\LRo. In generic terms lower bound (LB) and upper bound (UB), we define the duality gap as $\text{DG}=|\text{LB}-\text{UB}|/|\text{UB}|$. To assess the quality of $\hat{Z}$ obtained by the Lagrangian heuristic relative to $Z_{\text{G}}$ obtained by the Gurobi optimizer, we define the \textit{relative performance metric}, $\text{RP}=(\hat{Z}-Z_{\text{G}})/Z_{\text{G}}$, unrestricted in sign. If $\text{RP}>0$, the optimum obtained by the Lagrangian heuristic outperforms that of the Gurobi optimizer. If $\text{RP}<0$, the optimum obtained by the Gurobi optimizer outperforms that of the Lagrangian heuristic. If $\text{RP}=0$, the obtained optimums of both methods are the same.

\subsection{Computational Experiment Results} \label{sec:results}
Out of 180 RCRP-ARC test instances, we present detailed analyses for 36 instances. Table~\ref{tab:computational_results1} reports the computational results for test instances with $\varepsilon/\varepsilon_{\max}=0.3$, illustrating scenarios where resources are limited. For 9 ``small'' instances,  the baseline Gurobi optimizer successfully identified optimal solutions, or those within the specified duality gap tolerance of \SI{0.5}{\%}, within the specified runtime limit of \SI{3600}{s}. However, as the size of instances increases, we start to observe the Gurobi optimizer triggering the runtime limit. Particularly, for instances 17 and 18, we see a significant duality gap of \SI{11.03}{\%} and \SI{70.15}{\%}, respectively. Examining the results of the Lagrangian heuristic, we observe that all 18 instances were solved in less than \SI{462.24}{s}. Comparing the feasible primal solutions to RCRP-ARC, there are 10 instances in which Gurobi solutions performed better than the Lagrangian heuristic solutions. However, the differences are at most \SI{1.77}{\%}. The Lagrangian heuristic outperformed the Gurobi optimizer for 6 instances with the largest recorded margin of \SI{26.21}{\%} (instance 18) and obtained optimal solutions for 2 instances.

Table~\ref{tab:computational_results2} presents the results for cases where resources are abundant, specifically for instances with $\varepsilon/\varepsilon_{\max}=0.8$. All parameters are the same as those we have shown previously in Table~\ref{tab:computational_results1}, except for the $\varepsilon$ value. An increase in the value of $\varepsilon/\varepsilon_{\max}$ leads to an enlargement of the feasible solution set. Out of 18 instances, the Gurobi optimizer only solved one instance (instance 5) to the optimality within the runtime limit and one instance (instance 13) to the tolerance-optimality by the runtime limit. The Lagrangian heuristic solved all instances with the maximum runtime of \SI{810.97}{s}. The duality gaps obtained by the Lagrangian heuristic are comparably larger than those with the lower $\varepsilon$ because $Z_{\text{D}}$ converges to $Z_{\text{LP}}$. However, it is important to note that the Lagrangian relaxation bound is theoretically no worse than the LP relaxation bound (assuming converged multipliers), which suggests that the integrality gap of the problem is significant. For 12 out of 18 instances, the primal solutions obtained by the Lagrangian heuristic outperform the best incumbent primal solution of the Gurobi optimizer found by the runtime limit. The outperformance of the Lagrangian heuristic over the Gurobi optimizer is notably significant for instances~14--18, with the relative performance metric ranging from \SI{11.68}{\%} to \SI{25.84}{\%}. The underperformance of the Lagrangian heuristic is also observable, with the relative performance metric ranging up to \SI{1.36}{\%}.

\begin{table}[h]
	\fontsize{9}{10}\selectfont
	\renewcommand{\arraystretch}{1.2}
	\caption{Computational results for RCRP-ARC test instances with $\varepsilon/\varepsilon_{\max}=0.3$.}
	\centering
	\begin{threeparttable}
	\begin{adjustbox}{center}
	\begin{tabular}{@{\extracolsep{-5.8pt}}*{12}{r}@{}}
		\hline
		\hline
		\multirow{2}{*}{Instance} & \multicolumn{4}{c}{Gurobi (8-core)\tnote{*}} & \multicolumn{6}{c}{Lagrangian Heuristic} & \multirow{2}{*}{RP\tnote{\S}, \%} \\
		\cmidrule(lr){2-5} \cmidrule(lr){6-11}
		& LB & UB ($Z_{\text{G}}$) & Runtime\tnote{\textdagger}, s & DG, \% & LB ($Z_{\text{D}}$) & Runtime, s & UB ($\hat{Z}$) & Runtime, s & Total runtime\tnote{\textdaggerdbl}, s & DG, \% & \\
		\hline
		1 & -2,701 & -2,701 & 58.42 & 0 & -2,790.36 & 3.81 & -2,701 & 3.54 & 23.21 & 3.31 & 0 \\
2 & -2,860 & -2,855 & 20.04 & 0.18 & -2,908.58 & 25.06 & -2,839 & 10.98 & 57.29 & 2.45 & -0.56 \\
3 & -8,273 & -8,246 & 626.01 & 0.33 & -8,633.18 & 6.98 & -8,229 & 18.64 & 45.92 & 4.91 & -0.21 \\
4 & -7,759 & -7,721 & 1,765.20 & 0.49 & -8,026.34 & 24.36 & -7,658 & 25.27 & 71.21 & 4.81 & -0.82 \\
5 & -1,910 & -1,910 & 7.61 & 0 & -1,947.42 & 25.92 & -1,910 & 20.49 & 75.00 & 1.96 & 0 \\
\arrayrulecolor{gray!50}\hline
6 & -5,926 & -5,898 & 439.14 & 0.47 & -6,139.74 & 55.39 & -5,819 & 20.04 & 97.56 & 5.51 & -1.34 \\
7 & -10,414 & -10,372 & 1,039.13 & 0.40 & -10,952.97 & 6.82 & -10,357 & 24.57 & 51.89 & 5.75 & -0.14 \\
8 & -13,412 & -13,350 & 1,060.46 & 0.46 & -13,559.13 & 159.60 & -13,291 & 50.35 & 245.75 & 2.02 & -0.44 \\
9 & -18,530 & -17,704 & - & 4.67 & -18,565.75 & 22.50 & -17,976 & 98.79 & 155.67 & 3.28 & 1.54 \\
10 & -16,764 & -16,681 & 2,842.39 & 0.50 & -17,047.76 & 54.95 & -16,385 & 31.50 & 112.65 & 4.04 & -1.77 \\
\hline
11 & -15,549 & -15,456 & - & 0.60 & -15,719.22 & 48.23 & -15,436 & 70.60 & 149.94 & 1.83 & -0.13 \\
12 & -16,788 & -15,749 & - & 6.60 & -16,838.56 & 66.40 & -16,235 & 98.89 & 201.11 & 3.72 & 3.09 \\
13 & -13,570 & -13,376 & - & 1.45 & -13,641.58 & 14.10 & -13,328 & 51.32 & 87.62 & 2.35 & -0.36 \\
14 & -27,064 & -25,889 & - & 4.54 & -27,212.63 & 40.86 & -25,846 & 81.93 & 162.77 & 5.29 & -0.17 \\
15 & -21,399 & -20,643 & - & 3.66 & -21,523.70 & 25.99 & -20,801 & 30.80 & 88.18 & 3.47 & 0.77 \\
\hline
16 & -32,524 & -29,872 & - & 8.88 & -32,602.67 & 68.43 & -30,824 & 336.10 & 462.24 & 5.77 & 3.19 \\
17 & -47,164 & -42,480 & - & 11.03 & -47,302.04 & 40.05 & -46,219 & 148.74 & 242.97 & 2.34 & 8.80 \\
18 & -60,000 & -35,263 & - & 70.15 & -46,708.36 & 37.08 & -44,504 & 74.25 & 162.22 & 4.95 & 26.21 \\
		\arrayrulecolor{black}\hline
		\hline
	\end{tabular}
	\end{adjustbox}
	\begin{tablenotes}
        \item[*] Gurobi optimizer utilizes a combination of branch-and-bound, cutting planes, presolve, heuristics, and parallelism.
	\item[\textdagger] Hyphen (-) indicates the trigger of the runtime limit of \SI{3600}{s}.
	\item[\textdaggerdbl] Includes runtimes for LB, UB, and the subgradient method intermediate steps.
	\item[\S] Positive RP indicates the outperformance of the Lagrangian heuristic.
	\end{tablenotes}
	\end{threeparttable}
	\label{tab:computational_results1}
\end{table}
\begin{table}[h]
	\fontsize{9}{10}\selectfont
	\renewcommand{\arraystretch}{1.2}
	\caption{Computational results for RCRP-ARC test instances with $\varepsilon/\varepsilon_{\max}=0.8$.}
	\centering
	\begin{threeparttable}
	\begin{adjustbox}{center}
	\begin{tabular}{@{\extracolsep{-5.8pt}}*{12}{r}@{}}
		\hline
		\hline
		\multirow{2}{*}{Instance} & \multicolumn{4}{c}{Gurobi (8-core)\tnote{*}} & \multicolumn{6}{c}{Lagrangian Heuristic} & \multirow{2}{*}{RP\tnote{\S}, \%} \\
		\cmidrule(lr){2-5} \cmidrule(lr){6-11}
		& LB & UB ($Z_{\text{G}}$) & Runtime\tnote{\textdagger}, s & DG, \% & LB ($Z_{\text{D}}$) & Runtime, s & UB ($\hat{Z}$) & Runtime, s & Total runtime\tnote{\textdaggerdbl}, s & DG, \% & \\
		\hline
		1 & -2,800 & -2,750 & - & 1.82 & -2,794.30 & 13.01 & -2,725 & 12.24 & 39.15 & 2.54 & -0.91 \\
2 & -3,002 & -2,949 & - & 1.80 & -3,036.67 & 36.14 & -2,915 & 37.48 & 97.52 & 4.17 & -1.15 \\
3 & -8,998 & -8,404 & - & 7.07 & -9,079.56 & 23.18 & -8,368 & 43.12 & 85.05 & 8.50 & -0.43 \\
4 & -8,539 & -8,008 & - & 6.63 & -8,575.73 & 30.51 & -7,964 & 110.34 & 164.22 & 7.68 & -0.55 \\
5 & -1,910 & -1,910 & 24.22 & 0 & -1,917.86 & 32.89 & -1,910 & 22.45 & 74.05 & 0.41 & 0 \\
\arrayrulecolor{gray!50}\hline
6 & -7,105 & -6,324 & - & 12.35 & -7,139.11 & 72.30 & -6,332 & 65.31 & 163.22 & 12.75 & 0.13 \\
7 & -11,308 & -10,936 & - & 3.40 & -11,336.61 & 10.68 & -10,787 & 39.28 & 63.85 & 5.10 & -1.36 \\
8 & -14,595 & -13,888 & - & 5.09 & -14,853.03 & 39.79 & -13,937 & 120.63 & 184.69 & 6.57 & 0.35 \\
9 & -19,706 & -17,765 & - & 10.93 & -20,000.00 & 12.25 & -18,206 & 52.93 & 74.87 & 9.85 & 2.48 \\
10 & -19,561 & -16,667 & - & 17.36 & -19,855.75 & 42.07 & -17,499 & 105.23 & 167.91 & 13.47 & 4.99 \\
\hline
11 & -17,067 & -16,007 & - & 6.62 & -17,154.73 & 52.43 & -16,241 & 88.30 & 162.16 & 5.63 & 1.46 \\
12 & -17,722 & -15,963 & - & 11.02 & -17,804.39 & 138.79 & -16,819 & 274.64 & 457.87 & 5.86 & 5.36 \\
13 & -13,600 & -13,542 & - & 0.43 & -13,651.54 & 24.38 & -13,490 & 87.55 & 133.38 & 1.20 & -0.38 \\
14 & -30,000 & -21,791 & - & 37.67 & -29,339.01 & 65.11 & -27,033 & 316.57 & 427.93 & 8.53 & 24.06 \\
15 & -40,000 & -17,032 & - & 134.85 & -21,685.08 & 56.67 & -20,908 & 60.93 & 157.88 & 3.72 & 22.76 \\
\hline
16 & -40,000 & -28,719 & - & 39.28 & -34,860.21 & 113.75 & -32,073 & 643.39 & 810.97 & 8.69 & 11.68 \\
17 & -60,000 & -42,761 & - & 40.31 & -53,278.35 & 72.80 & -48,932 & 399.13 & 550.60 & 8.88 & 14.43 \\
18 & -60,000 & -36,585 & - & 64.00 & -49,433.34 & 73.83 & -46,037 & 232.74 & 366.34 & 7.38 & 25.84 \\
		\arrayrulecolor{black}\hline
		\hline
	\end{tabular}
	\end{adjustbox}
	\begin{tablenotes}
        \item[*] Gurobi optimizer utilizes a combination of branch-and-bound, cutting planes, presolve, heuristics, and parallelism.
	\item[\textdagger] Hyphen (-) indicates the trigger of the runtime limit of \SI{3600}{s}.
	\item[\textdaggerdbl] Includes runtimes for LB, UB, and the subgradient method intermediate steps.
	\item[\S] Positive RP indicates the outperformance of the Lagrangian heuristic.
	\end{tablenotes}
	\end{threeparttable}
	\label{tab:computational_results2}
\end{table}

\begin{figure}[htbp]
	\centering
		\includegraphics[width=0.68\linewidth]{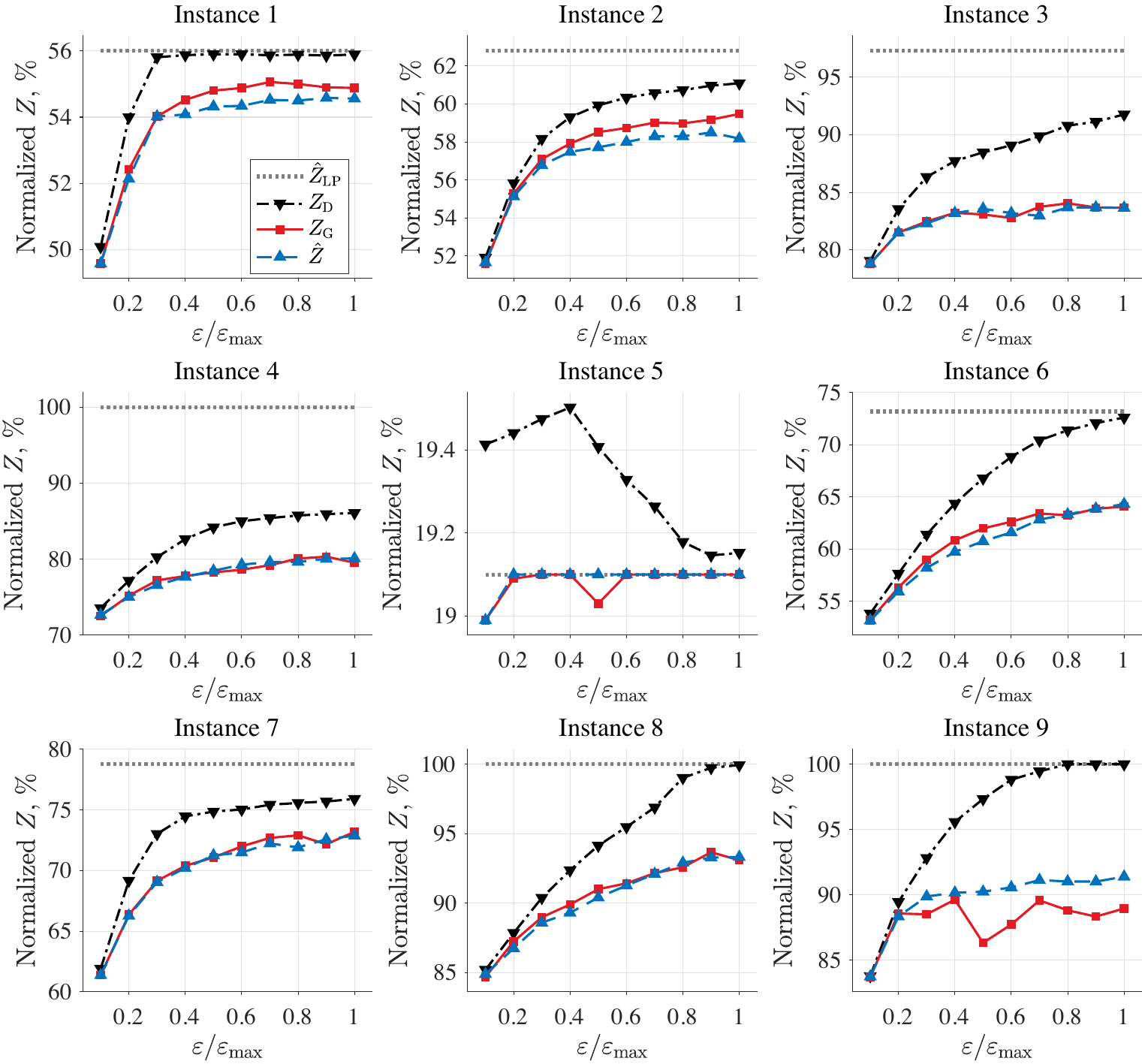}
		\caption{Computational results for instances 1--9. Note that all metrics are normalized and flipped in sign.}
		\label{fig:pareto1}
\end{figure}

\begin{figure}[htbp]
	\centering
		\includegraphics[width=0.68\linewidth]{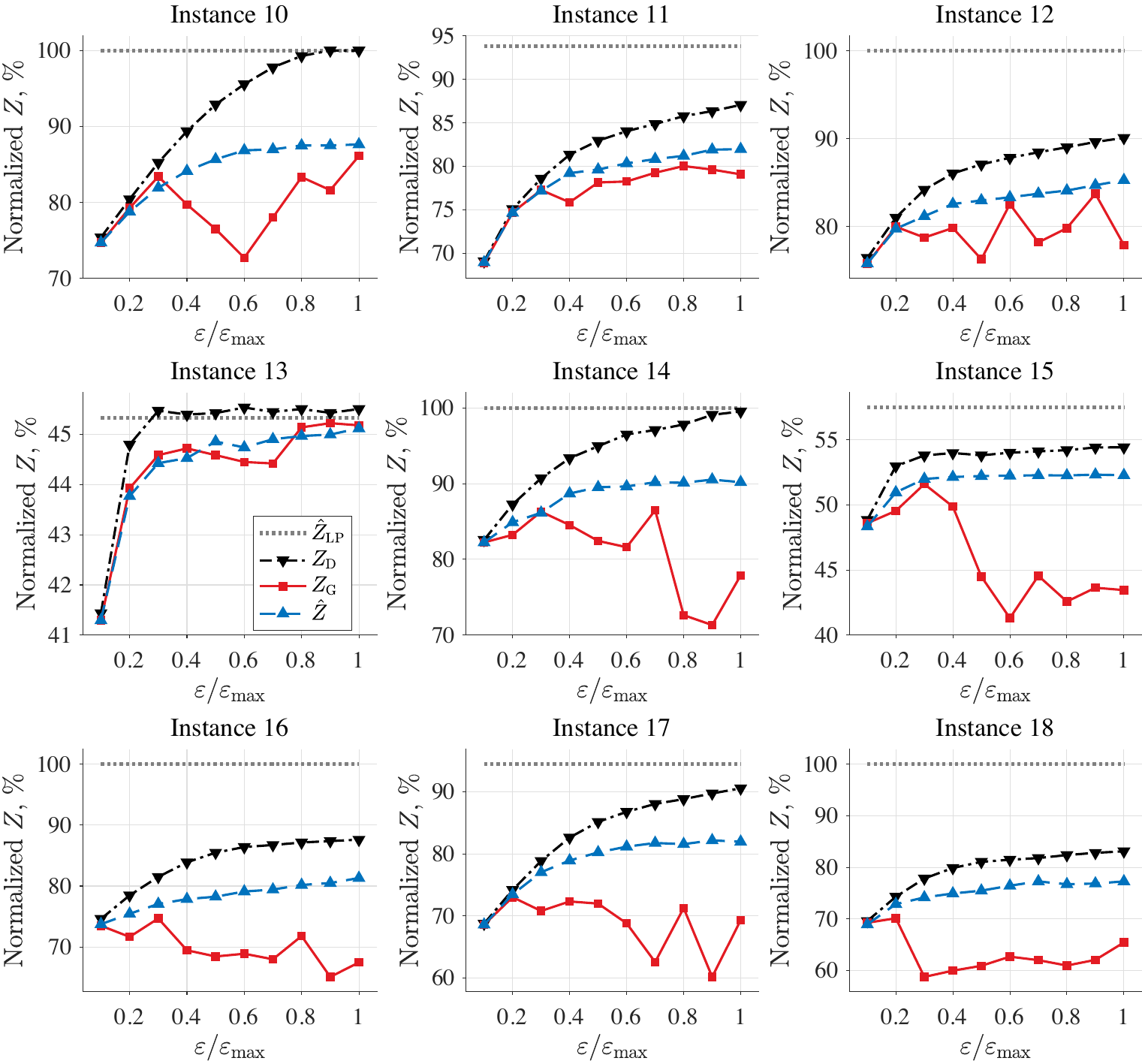}
		\caption{Computational results for instances 10--18. Note that all metrics are normalized and flipped in sign.}
		\label{fig:pareto2}
\end{figure}

\begin{figure}[htbp]
	\centering
	\includegraphics[width=0.68\linewidth]{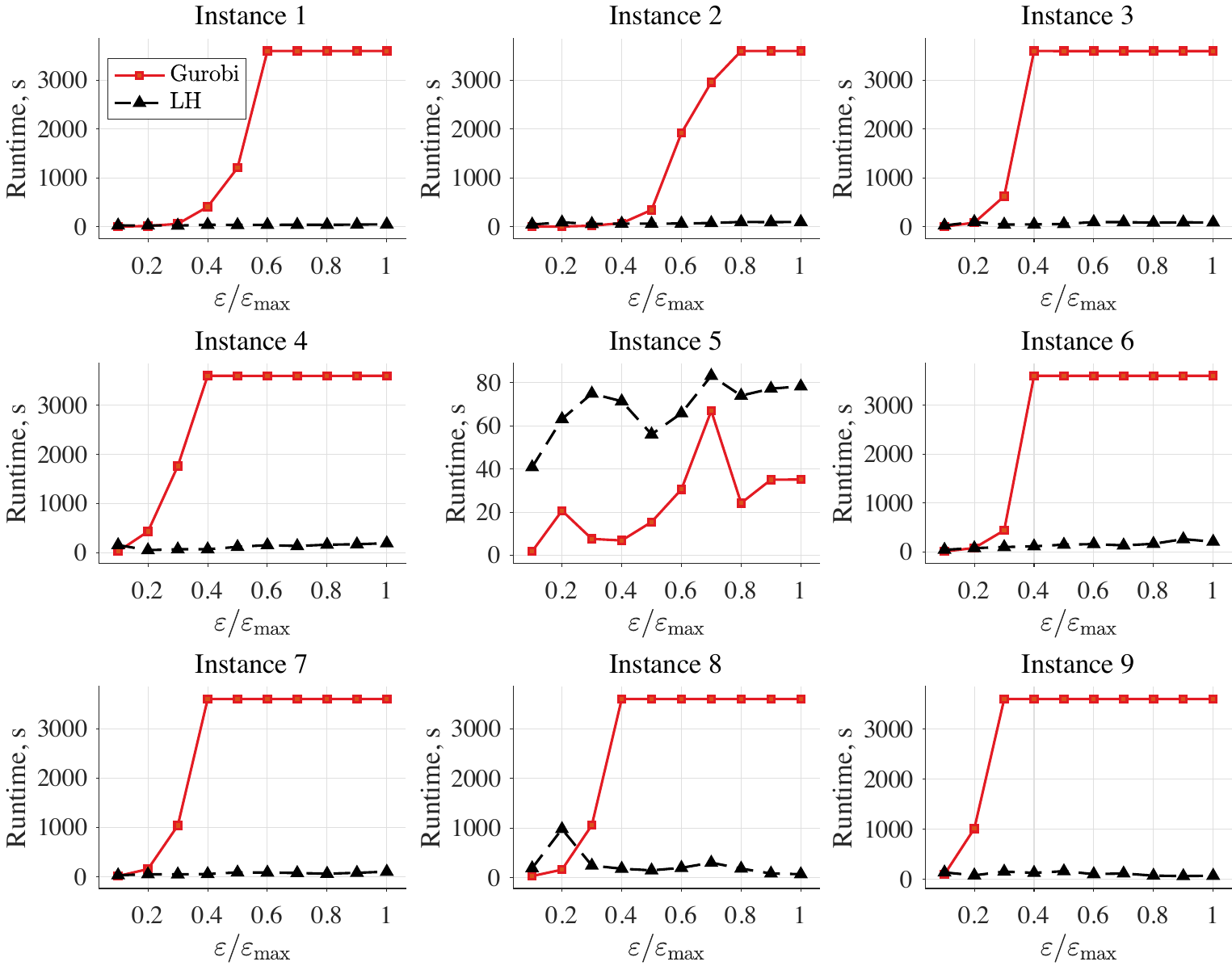}
	\caption{Runtime results for instances 1--9. The runtime limit of \SI{3600}{s} is enforced.}
	\label{fig:time1}
\end{figure}

\begin{figure}[htbp]
	\centering
		\includegraphics[width=0.68\linewidth]{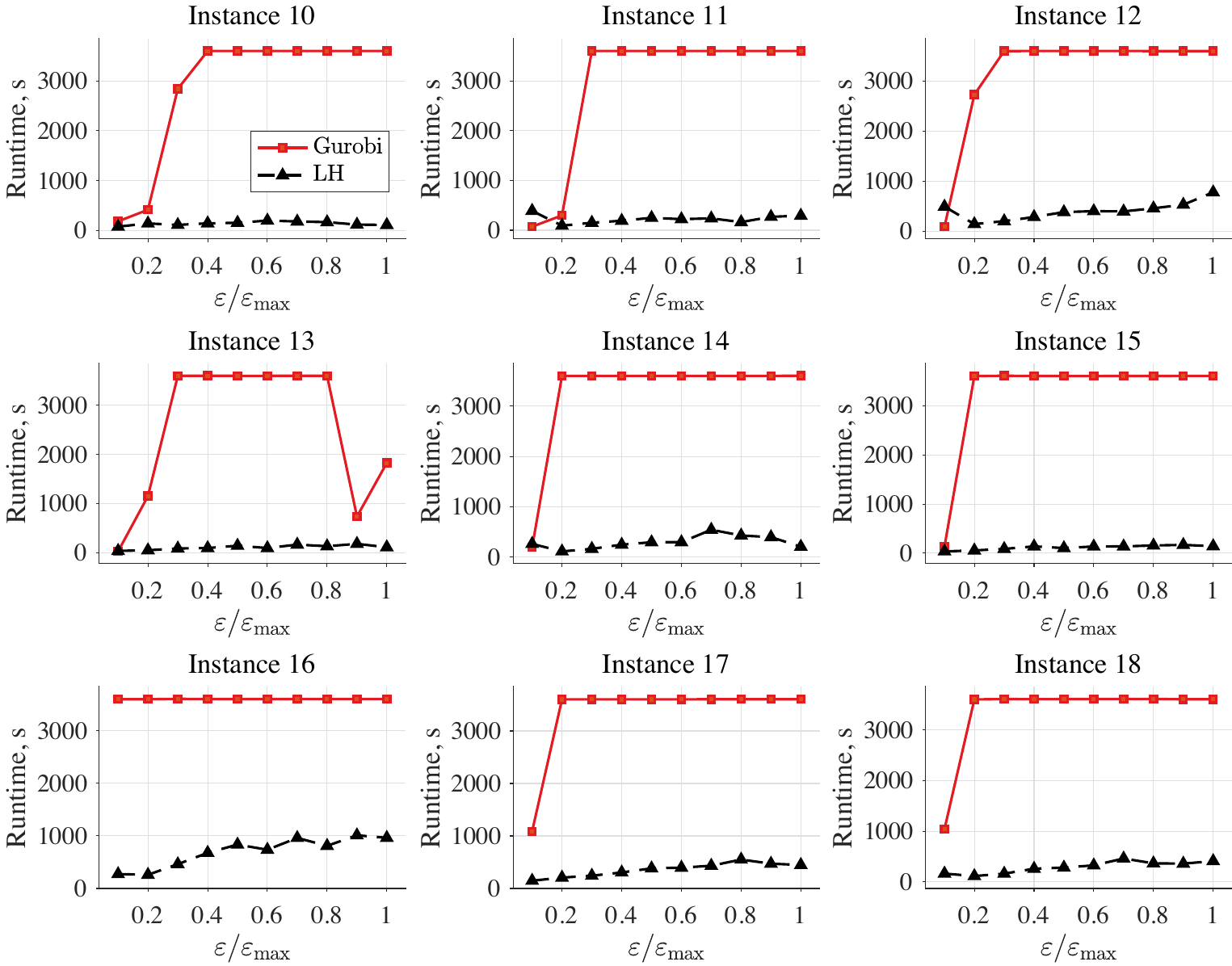}
		\caption{Runtime results for instances 10--18. The runtime limit of \SI{3600}{s} is enforced.}
		\label{fig:time2}
\end{figure}

We present the results of all 180 test instances graphically. Figures~\ref{fig:pareto1} and \ref{fig:pareto2} visualize the computational results, showcasing the approximated Pareto fronts of both methods. In these figures, all metrics, $\hat{Z}_{\text{LP}}$, $Z_{\text{D}}$, $Z_{\text{G}}$, and $\hat{Z}$, are normalized and have their signs flipped for ease of physical interpretation. Note that generating the true Pareto front of a given RCRP instance requires solving all associated RCRP-ARC instances to optimality.\footnote{Strictly speaking, the Pareto front in the discrete-time domain is also the approximation of the true Pareto front in the continuous-time domain.} Without an optimality certificate (which is typically proven by the duality gap), the results, $Z_{\text{D}}$ and $Z_{\text{G}}$, in Figs.~\ref{fig:pareto1} and \ref{fig:pareto2} are deemed approximations of Pareto fronts; the dominated solutions are still included for completeness. Figure~\ref{fig:pareto2} corroborates the outperformance of the Lagrangian heuristic for large instances. For instances 1, 5, and 13, we observe that $\hat{Z}_{\text{LP}}$ effectively certifies that $Z_{\text{D}}$ has either converged or is not optimal, albeit its usefulness appears to be limited. Figures~\ref{fig:time1} and \ref{fig:time2} compare the computation runtime between the 8-core Gurobi optimizer and the Lagrangian heuristic with no parallel computing implementation (except that we solve \LRo using the 8-core Gurobi optimizer). In most cases, we see that the Gurobi optimizer reached the runtime limit of \SI{3600}{s}. A notable case is instance~5, in which the initial solution is near-optimal and no significant maneuvers are needed to maximize the total coverage reward.\footnote{The Gurobi optimizer returned tolerance-optimal (\SI{0.5}{\%} solutions for instance 5, hence, we observe a drop in $\varepsilon/\varepsilon_{\max}=0.5$ RCRP-ARC.}



\subsection{Illustrative Example: Federated Disaster Monitoring} \label{sec:illustrative_example}
In this illustrative example, we focus on demonstrating the versatility of the RCRP formulation and solution method by showcasing (i) the application of the RCRP-IRC formulation, (ii) the utilization of the Lagrangian heuristic method in conjunction with RCRP-IRC, and (iii) non-RGT orbits.

\subsubsection{Problem Setup}
Suppose a group of seven satellites in different circular orbits (parameters shown in the left half of Table~\ref{tab:rcrp-irc}) is tasked with a reconfiguration process to form a federation for a 15-day satellite-based emergency mapping mission to monitor active disaster events and support post-disaster relief operations. The spot targets of interest are Getty, California $(34.09\degree \text{N},118.47\degree \text{W})$, Asheikri, Nigeria $(11.96\degree \text{N},12.93\degree \text{E})$, and Hunga Tonga–Hunga Ha'apai, Tonga $(21.18\degree \text{S},175.19\degree \text{W})$. The coverage rewards are randomly generated following the standard uniform distribution in the range of $[0,1]$. We also let $r_{tp}=1,\forall t\in\mathcal{T},p\in\mathcal{P}$ and all targets enforce $\vartheta_{\min}=10\degree$.

We consider a set of orbital slots, $\mathcal{J}=\{\mathcal{J}_1,\ldots,\mathcal{J}_7\}$, where $\mathcal{J}_i$ denotes the set of orbital slots that are $\Delta v$-compatible with satellite $i$. This means that the cost of transferring satellite $i$ to orbital slot $j$ in $\mathcal{J}_i$ is less than or equal to the specified $\Delta v$ value. Each $\mathcal{J}_i$ comprises orbital slots that allow satellite $i$ to perform one of the following four options: (i) change inclination, (ii) change RAAN, (iii) make a coplanar phasing maneuver, or (iv) stay in its orbit. To generate orbital slots for the first option, we determine the boundary inclination values and generate inclinations uniformly within the range of $[inc_{i,\text{LB}},inc_{i,\text{UB}}]$ given $\varepsilon_i$. Similarly, for the second option, we find the boundary values and generate RAAN values uniformly within the range of $[\Omega_{i,\text{LB}},\Omega_{i,\text{UB}}]$. For the phasing maneuver, orbital elements of orbital slots are the same as the satellites except for the arguments of latitude $u$, which is uniformly distributed in $[0,360\degree)$. Lastly, we add $|\mathcal{I}|$ initial orbits to $\mathcal{J}$ to allow for no maneuvering option for satellites. Globally, regardless of maneuvering options, we let $a_j=a_i$ and $e_j=0$. For the first and second options, we generate $u_j$ uniformly distributed in the range of $[0,360\degree)$ similar to the phasing maneuver option. Note that while $\mathcal{J}_i$ is constructed for satellite $i$, other satellites beside satellite $i$ may transfer to orbital slots in $\mathcal{J}_i$ as long as the resource constraints are not violated.

For orbit propagation, we use the SGP4 (Simplified General Perturbations-4) model to account for the differential secular and periodic rates in the change of orbital elements that affect satellite states and visibility matrix during the entire specified mission time horizon $T=\SI{15}{days}$. It is important to note that satellites will continue to be subject to differential orbit perturbations beyond the considered time horizon, and thus, station-keeping maneuvers may be warranted if repeatability of the coverage state is desired. For the cost matrix and to account for the possibility of an altitude change, we use the combined plane change and the Hohmann transfer maneuvers and the coplanar phasing maneuvers as outlined in Chapters 6.5.1 and 6.6.1 in Ref.~\cite{vallado2013fundamentals}, respectively. The phasing angle is set to $180\degree$ as the worst-case value, and thus, the corresponding $\Delta v_{\text{phasing}}^\prime$ serves as the upper limit for the actual $\Delta v_{\text{phasing}}$ required for phasing for a given phasing time. This assumption is made to account for potential inaccuracies in phasing modeling. This allows us to decouple plane changes (inclination and RAAN changes) from phasing, and hence, we only need to enforce resource constraints for plane change maneuvers in the optimization. To calculate the complete $\Delta v$, we add $\Delta v_{\text{phasing}}^\prime$ to $\Delta v_{\text{pc}}$, the cost of the plane change maneuver. It is important to factor in the phasing cost when determining the $\varepsilon$ allocated for a plane change.

The size of the RCRP-IRC instance is as follows. We let $|\mathcal{J}_i|=1,801,\forall i \in\mathcal{I}$, and thus, we have $|\mathcal{J}|=12,607$. We also let the time step size to be \SI{120}{s}, and consequently, $|\mathcal{T}|=10,800$. With $|\mathcal{I}|=7$ and $|\mathcal{P}|=3$, the instance has 88,249 assignment variables, 32,400 coverage state variables, and 45,021 constraints (excluding the decision variable domain definitions).

\subsubsection{Numerical Results}
Letting $\varepsilon_i=\SI{1}{km/s},\forall i \in\mathcal{I}$ as the $\Delta v$ budgets allocated for plane change maneuvers, and using the Lagrangian heuristic method with a neighborhood size of $|\mathcal{N}^\prime|\le 50|\mathcal{I}|$ to solve RCRP-IRC, we obtained $\hat{Z}=-3198.19$ in \SI{329.60}{s} with the duality gap of \SI{8.01}{\%}. Note that the initial configuration has a score of $Z=-2632.75$. The LH optimum improved the initial score by \SI{21.48}{\%}. The final configuration is specified in Table~\ref{tab:rcrp-irc} (right half). The $\Delta v_{\text{pc}}$ column shows the results from the plane change maneuvers\footnote{To assess the overall $\Delta v$, we would need to add $\Delta v_{\text{phasing}}$, which can be obtained by trading off the time required for the phasing (the longer the phasing time, the lower the phasing cost, and vice versa).}. 

\begin{figure}[htbp]
	\centering
		\includegraphics[width=0.5\linewidth]{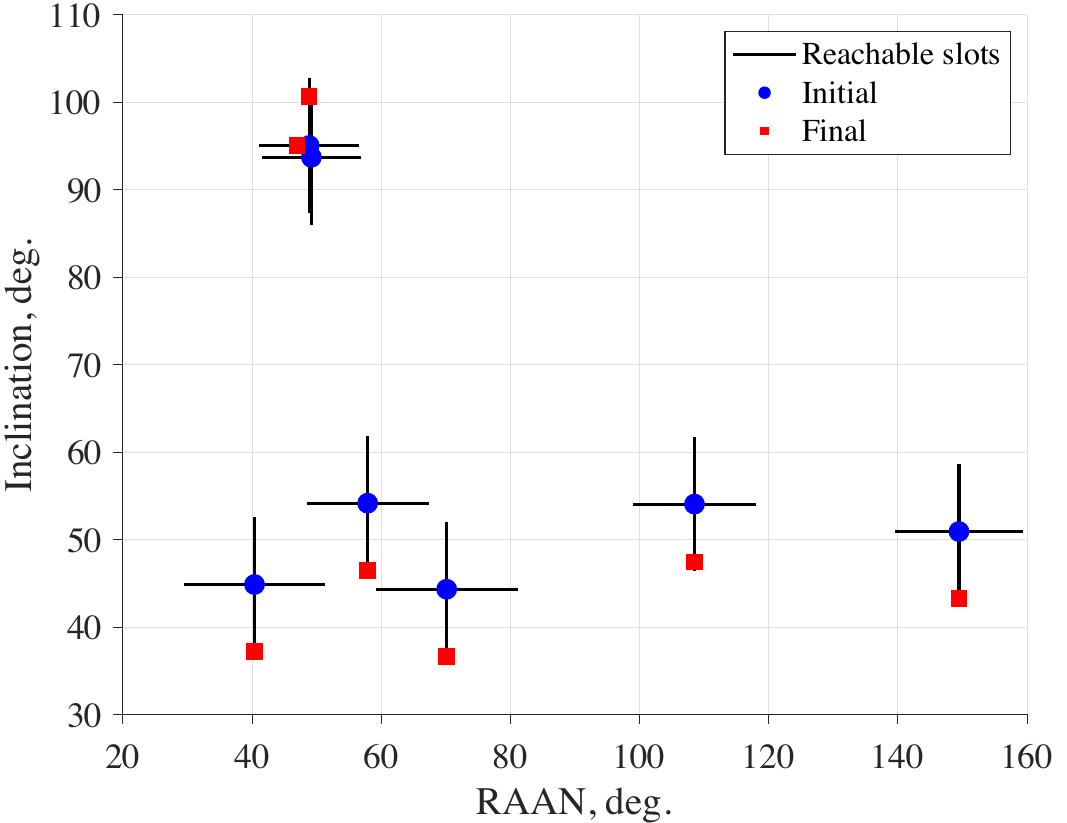}
		\caption{Range of reachable orbital slots (shown in line), initial (blue circle) and final (red square) configurations.}
		\label{fig:incvsraan}
\end{figure}

Interestingly, not all satellites fully used their allocated $\Delta v_{\text{pc}}$ budgets even though there was no penalty for using it up to the limit. Satellites 2--6 lowered their inclinations to maximize coverage, as the targets were in a low-latitude zone. Satellites 2--5 performed the maximum inclination change possible. In contrast, satellite 1 only changed its RAAN, despite having a near-polar inclination. Satellite 7 was assigned to an orbital slot generated for satellite 1; the maneuver was possible due to the close proximity in $\Delta v$ required. This resulted in a decrease in altitude by \SI{14}{km} and a slight change in RAAN, where the former effectively reduced the swath width of a sensor. Figure~\ref{fig:incvsraan} illustrates the reconfiguration process results in $(\Omega,inc)$ space. The horizontal and vertical lines each respectively indicate the range of RAAN and inclination values a satellite (in blue) is allowed to transfer given the resource constraint.

\begin{figure}[htbp]
	\centering
		\includegraphics[width=0.5\linewidth]{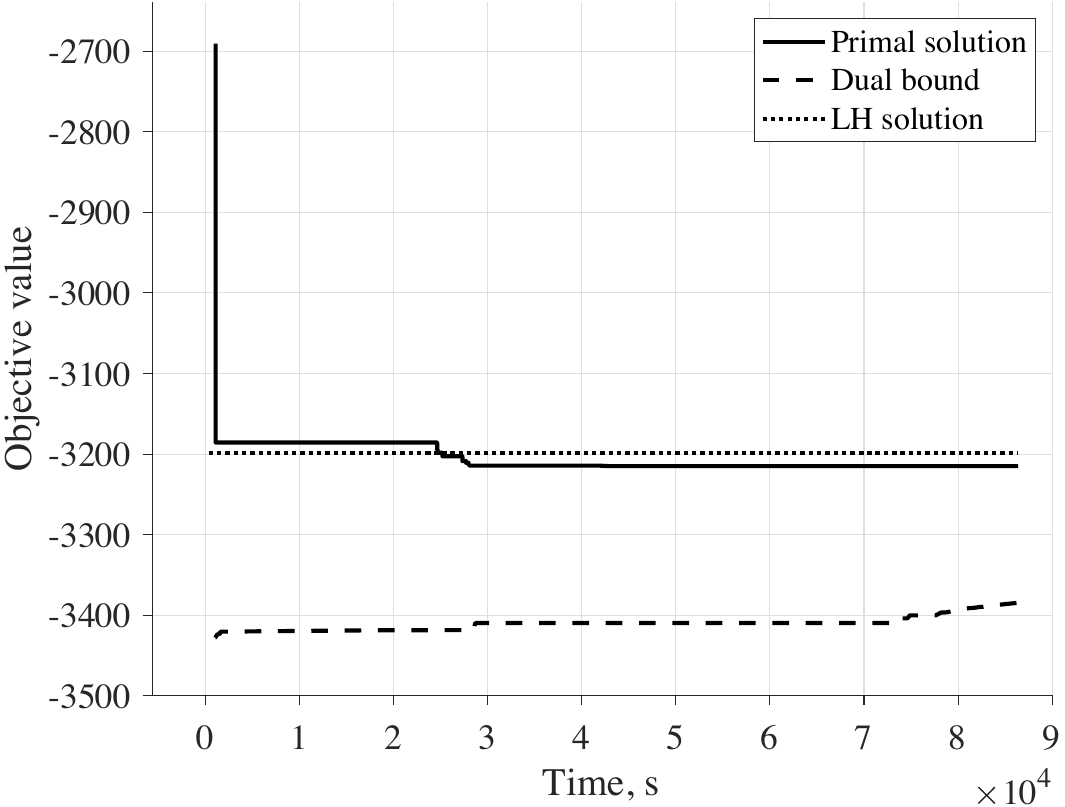}
		\caption{Primal solution and dual bound over time.}
		\label{fig:time}
\end{figure}

By solving the same instance of RCRP-IRC with the Gurobi optimizer using the same setting as in Section~\ref{sec:setup}, we obtained $Z_{\text{G}}=-3185.51$ at the end of the runtime limit of \SI{3600}{s}, with a duality gap of \SI{7.32}{\%}. This solution is underperforming compared to $\hat{Z}$. When the runtime limit was extended to \SI{86400}{s}, the Gurobi optimizer resulted in $Z_{\text{G}}=-3214.934$ with a duality gap of \SI{5.26}{\%}. This solution improved the LH solution by \SI{0.52}{\%}. As the Gurobi optimizer was given more time to converge, it found a more optimal solution than the LH method after \SI{24660}{s} into the optimization. However, the large scale of the problem prevented the Gurobi optimizer from finding the optimal solution or proving the optimality of the incumbent solution by closing the gap within a day of runtime. The progress of the Gurobi optimizer is shown in Fig.~\ref{fig:time}.

Finally, we conducted an additional experiment with only phasing maneuvers. The both LH and Gurobi optimizer produced the same solution, $\hat{Z}=-2849.78$, in \SI{70.18}{s} and \SI{1765.68}{s}, respectively. The Gurobi optimizer spent extra time proving the optimality of the solution. The obtained solution improved the initial configuration by \SI{8.24}{\%}; this result corroborates the effectiveness of reconfiguration, even when it involves performing only rephasing among satellites. The obtained optimal solution is as follows: $u_1=48\degree,u_2=192\degree,u_3=60\degree,u_4=96\degree,u_5=168\degree,u_6=336\degree$, and $u_7=348\degree$.

\begin{table}[htbp]
    \fontsize{9}{10}\selectfont
    \renewcommand{\arraystretch}{1.2}
	\caption{Problem setting (left) and obtained Lagrangian heuristic solution (right). The values in the parenthesis indicate the change in value and all orbital elements are defined at the epoch, J2000.}
	\centering
	\begin{tabular}{@{}r r r r r r r r r r@{}}
	\hline
	\hline
	\multirow{2}{*}{Satellite} & \multicolumn{4}{c}{Initial configuration}& \multicolumn{4}{c}{Final configuration (LH solution)} & \multirow{2}{*}{$\Delta v_{\text{pc}}$, km/s} \\
	\cmidrule(lr){2-5} \cmidrule(lr){6-9}
	& $a$, km & $inc$, deg. & $\Omega$, deg. & $u$, deg. & $a$, km & $inc$, deg. & $\Omega$, deg. & $u$, deg. & \\
	\hline
	1 & 7,161.83 & 95.04 & 48.83  & 275.54 & 7,161.83         & 95.04          & 46.96 (-1.86) & 48.00  & 0.24 \\
2 & 7,175.16 & 44.33 & 70.13  & 216.91 & 7,175.16         & 36.63 (-7.69)  & 70.13         & 0.00   & 1.00 \\
3 & 7,122.86 & 54.16 & 57.89  & 222.07 & 7,122.86         & 46.49 (-7.66)  & 57.89         & 192.00 & 1.00 \\
4 & 7,106.60 & 44.87 & 40.38  & 240.04 & 7,106.60         & 37.21 (-7.60)   & 40.38         & 228.00 & 1.00 \\
5 & 7,155.26 & 50.92 & 149.45 & 266.71 & 7,155.26         & 43.23 (-7.68)  & 149.45        & 324.00 & 1.00 \\
6 & 7,108.03 & 54.06 & 108.50 & 146.05 & 7,108.03         & 47.45 (-6.60)  & 108.50        & 336.00 & 0.86 \\
7 & 7,175.77 & 93.70 & 49.18  & 62.90  & 7161.83 (-13.94) & 100.60 (+6.91) & 48.83 (-0.35) & 204.00 & 0.90 \\
	\hline
	\hline
	\end{tabular}
	\label{tab:rcrp-irc}
\end{table}

\section{Conclusions} \label{sec:conclusion}
This paper proposes an integrated constellation design and transfer model for solving the RCRP. Given a set of target points each associated with a time-varying coverage reward and a time-varying coverage threshold, the problem aims to maximize the total reward obtained during a specified time horizon and to minimize the total cost of satellite transfers. The bi-objective formulation results in a trade-off analysis, potentially a Pareto front analysis if all $\varepsilon$ instances are solved to optimality, in the objective space spanned by the aggregated cost and the total coverage reward. Furthermore, as demonstrated in the illustrative example, the formulation can accommodate different types of orbits, not necessarily restricting orbital slots to be RGT orbits. The use of non-RGT orbits requires a user to specify the time horizon $T$ for which the formulation is valid.

The ILP formulation of RCRP-ARC enables users to utilize commercial software packages for convenient handling and obtaining tolerance-optimal solutions. However, for large-scale real-world instances, the problem suffers from the explosion of a combinatorial solution space. To overcome this challenge and to produce high-quality feasible primal solutions, we developed a Lagrangian relaxation-based heuristic method that combines the subgradient method with the 1-exchange neighborhood local search, exploiting the special substructure of the problem. The computational experiments in Section~\ref{sec:experiments} demonstrate the effectiveness of the proposed method, particularly for large-scale instances, producing near-optimal solutions with significantly reduced computational runtime compared to the reference solver.

We believe that the developed method provides an important step toward the realization of the concept of reconfiguration as a means for system adaptability and responsiveness, adding a new dimension to the operation of the next-generation satellite constellation systems.

\section*{Appendix A: Algorithms}
This appendix section lists the pseudocode of the algorithms discussed in the paper.
\begin{algorithm}
\DontPrintSemicolon
\caption{$\varepsilon$-constraint method} \label{alg:overall}
\KwIn{$\bm{c},\bm{\pi},\bm{r},\bm{v}$}
\KwOut{List of Z($\varepsilon$) values}
Initialize $\varepsilon \leftarrow \varepsilon_0$ by solving AP\;
\Repeat{termination flag is triggered}{
Z($\varepsilon$) $\leftarrow$ RCRP-ARC\;
$\varepsilon \leftarrow \varepsilon + \varepsilon_{\text{step}}$\;
}
\end{algorithm}

\begin{algorithm}
\DontPrintSemicolon
\caption{Subgradient optimization}
\KwIn{$\varepsilon,\bm{c},\bm{\pi},\bm{r},\bm{v}$}
\KwOut{$Z_{\text{D}}(\varepsilon),\hat{Z}(\varepsilon),(\bm{\varphi}^\ast$, $\tilde{\bm{y}}^\ast)$} \label{alg:subgradient}
$k \leftarrow 0$ and initialize $\bm{\lambda}^0$\;
\Repeat{termination flag is triggered}{
Solve $\text{LR}^k$: compute $Z_{\text{D}}(\varepsilon,\bm{\lambda}^k)$ and obtain optimal solution ($\bm{\varphi}^k,\bm{y}^k$) \Comment*[r]{LB; See Section~\ref{sec:lb}}
Compute $\hat{Z}(\varepsilon)$ and $(\bm{\varphi}^\ast,\tilde{\bm{y}}^\ast)$ via local search in the neighborhood $\mathcal{N}(\bm{\varphi}^k)$ of $\bm{\varphi}^k$ \Comment*[r]{UB; See Algorithm~\ref{alg:heuristic}}
Compute the subgradient $\bm{g}^k$ of $Z_{\text{D}}(\varepsilon,\bm{\lambda}^k)$ at $\bm{\lambda}^k$\;
Compute the step size $\theta_k$\;
Update Lagrange multipliers: $\bm{\lambda}^{k+1} \leftarrow \max(\bm{0},\bm{\lambda}^k+\theta_k\bm{g}^k)$\;
$k \leftarrow k+1$\;
}
\end{algorithm}

\begin{algorithm}
\DontPrintSemicolon
\caption{1-exchange neighborhood local search (full-scale)} \label{alg:heuristic}
\KwIn{$\varepsilon,\bm{c},\bm{\pi},\bm{r},\bm{v},\bm{\lambda},\bm{\varphi}$}
\KwOut{$\hat{Z}(\varepsilon),(\bm{\varphi}^\ast$, $\tilde{\bm{y}}^\ast(\bm{\varphi}^\ast))$}
Initialize $\mathcal{N}(\bm{\varphi})$\;
\Repeat{termination flag is triggered}{
Compute $\hat{Z}(\varepsilon)$ and $\bm{\varphi}^\ast$ from $\mathcal{N}(\bm{\varphi})$\;
Update $\bm{\varphi} \leftarrow \bm{\varphi}^\ast$\;
Update the neighborhood $\mathcal{N}(\bm{\varphi})$\;
}
Compute $\tilde{\bm{y}}^\ast(\bm{\varphi}^\ast)$ from $\bm{\varphi}^\ast$\;
\end{algorithm}

\section*{Appendix B: Selecting between Competing Relaxations} \label{sec:competing}
There exist different types of Lagrangian relaxations for RCRP, and the choice of $\varepsilon$-constraint transformation affects the complexity of the downstream algorithmic efforts and the mathematical properties. At first glance, one may observe that the coverage reward maximization objective function can be recast as an $\varepsilon$-constraint. Similar to the one proposed in this paper, the Lagrangian relaxation problem would be separable into two subproblems based on the type of variables. In such a case, the $\bm{\varphi}$ subproblem can be solved as an LP, and the $\bm{y}$ subproblem would be a relatively easy constrained ILP. Therefore, the lower bound calculation would still be computationally efficient. However, the main difference lies in the computation of $\hat{Z}(\varepsilon)$. Unlike the one discussed earlier, $\tilde{\bm{y}}^k(\bm{\varphi}^k)$ computed from $\bm{\varphi}^k$ would not necessarily satisfy the $\varepsilon$-constraint. Hence, additional considerations must come into play in obtaining a feasible primal solution. One viable approach is to solve the reduced formulation of RCRP, which fixes and parameterizes a subset of assignments from $\bm{\varphi}^k$ while optimizing the complement set. However, this approach becomes computationally expensive for instances with high $\varepsilon$ values. This approach was explored in our preliminary work \cite{lee2021lagrangian}.

One could attempt to relax an alternative set of constraints that may yield a tighter Lagrangian relaxation bound than the one proposed in this subsection. However, Lagrangian relaxation problems with Constraints~\eqref{eq:c3} retained may be unsuitable for embedding into an algorithm of iterative nature due to computational complexity. The rationale for the relaxation of Constraints~\eqref{eq:c3} can also be found in a study by Galvao and ReVelle \cite{galvao1996}, which reported a successful application of the Lagrangian relaxation of the linking constraints for MCLP; in their problem context, the Lagrangian relaxation problem possesses the integrality property.

\section*{Appendix C: Small RCRP Instance}
In this appendix section, we compare the performance of the developed Lagrangian heuristic method with the Gurobi optimizer, a mixed-integer programming solver. In this case, we set $|\mathcal{I}|=5$, $|\mathcal{J}|=200$, $|\mathcal{T}|=200$, and $|\mathcal{P}|=10$. Target points are randomly distributed within the maximum latitude bounds determined by the inclination of the satellites' orbits. Let $\textbf{\oe}_0=(a,e,inc,\Omega,u)=(\SI{8176.5}{km},0,75\degree,50\degree,0\degree)$. We also set $r_{tp}=1, \forall t\in\mathcal{T},p\in\mathcal{P}$ and all targets enforce $\vartheta_{\min}=7\degree$. $\mathcal{J}$ follows the common RGT constellation distribution rule. We use the default duality gap of \SI{0.01}{\%} for the Gurobi optimizer.

\begin{figure}[htbp]
	\centering
	\begin{subfigure}[h]{0.4\linewidth}
		\centering
		\includegraphics[width=\linewidth]{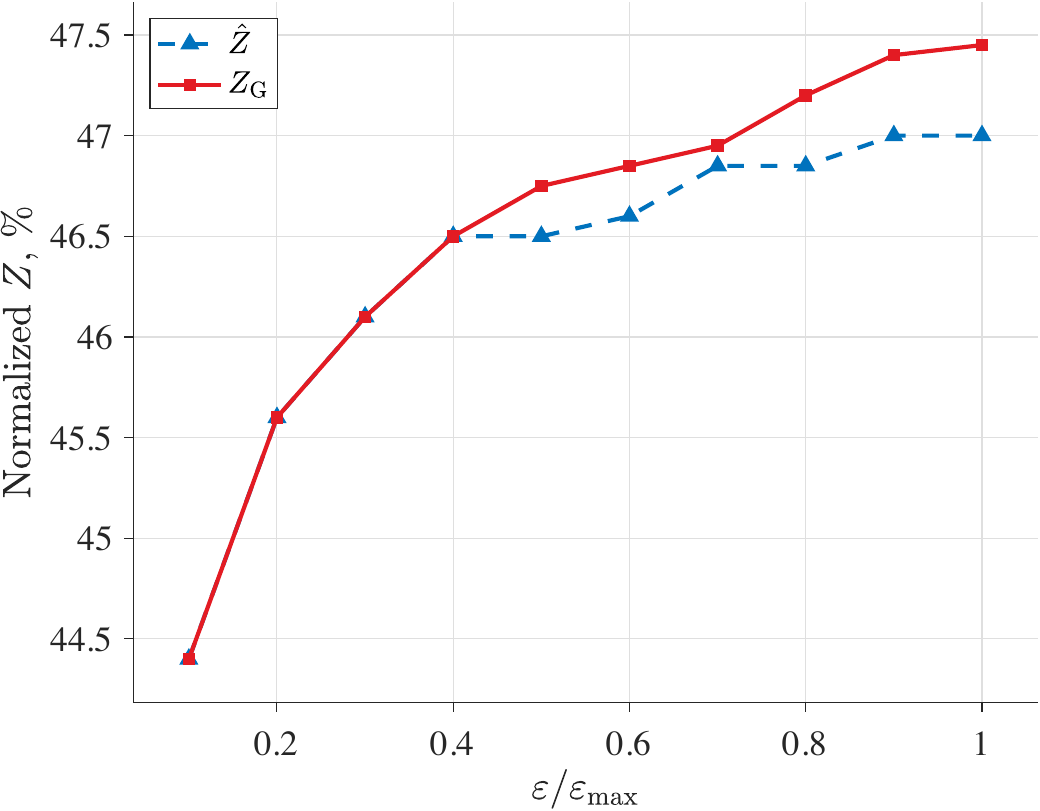}
		\caption{Pareto front comparison.}
	\end{subfigure}
	\begin{subfigure}[h]{0.4\linewidth}
		\centering
		\includegraphics[width=\linewidth]{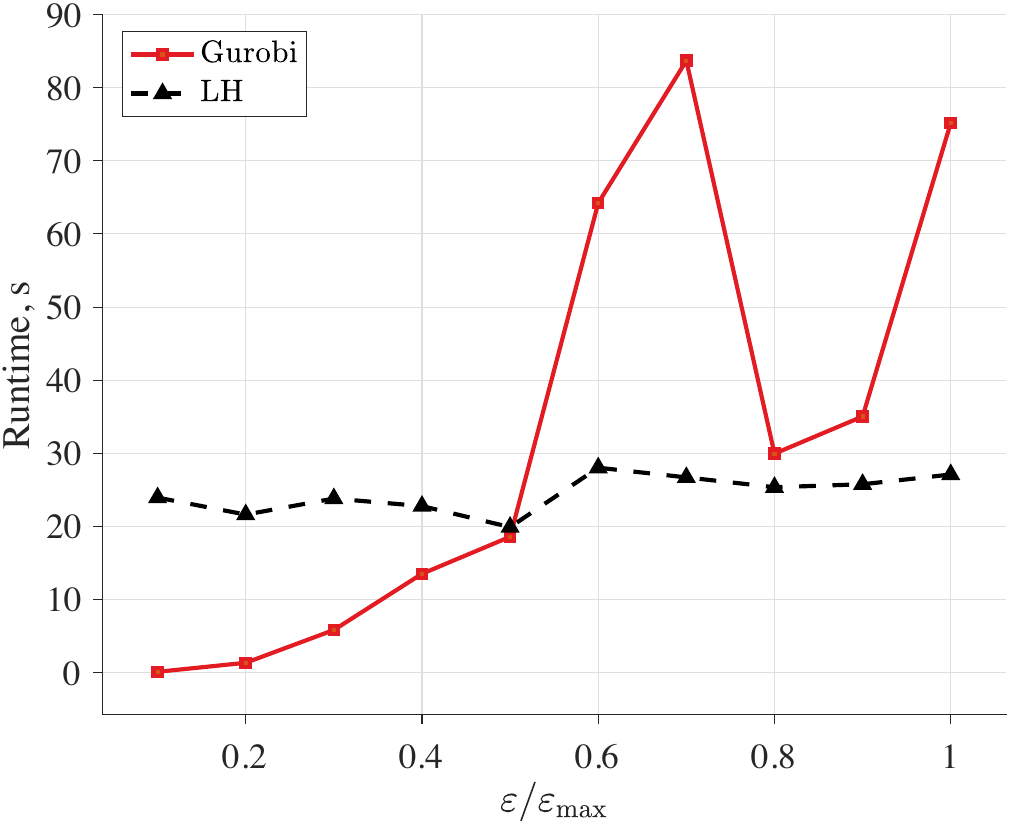}
		\caption{Runtime comparison.}
	\end{subfigure}
	\caption{Comparison of the Lagrangian heuristic method and the commercial optimizer for a small RCRP instance.}
	\label{fig:small}
\end{figure}

The results indicate that the Gurobi optimizer successfully converges to optimal solutions for all ten RCRP-ARC instances within the runtime limit, thereby identifying the true Pareto front of the RCRP. The maximum runtime reported is \SI{83.69}{s} for $\varepsilon/\varepsilon_{\max}=0.7$. The Lagrangian heuristic method does not find optimal solutions for RCRP-ARC with $\varepsilon/\varepsilon_{\max}\ge0.5$, and the maximum relative underperformance is \SI{0.95}{\%}. Pareto front and runtime analyses are shown in Fig.~\ref{fig:small}. These results demonstrate that, for small-scale instances, the use of conventional MILP methods can effectively characterize the Pareto front.

\section*{Acknowledgment}
This material is partially based upon work supported by the National Science Foundation Graduate Research Fellowship Program under Grant No. DGE-2039655. Any opinions, findings, and conclusions or recommendations expressed in this material are those of the author(s) and do not necessarily reflect the views of the National Science Foundation. The authors would like to express their gratitude to Sungwoo Kim and Da Eun Shim at Georgia Tech, and the anonymous reviewers, for their insightful suggestions.

\bibliography{references}

\begin{thebibliography}{49}
\newcommand{\enquote}[1]{``#1''}
\providecommand{\natexlab}[1]{#1}
\providecommand{\url}[1]{\texttt{#1}}
\providecommand{\urlprefix}{URL }
\expandafter\ifx\csname urlstyle\endcsname\relax
  \providecommand{\doi}[1]{\discretionary{}{}{}https://doi.org/#1}\else
  \providecommand{\doi}[1]{\discretionary{}{}{}\urlstyle{rm}\url{https://doi.org/#1}}\fi

\bibitem[{He et~al.(2020)He, Li, Yang, and Zhao}]{he2020reconfigurable}
He, X., Li, H., Yang, L., and Zhao, J., \enquote{{Reconfigurable Satellite
  Constellation Design for Disaster Monitoring Using Physical Programming},}
  \emph{International Journal of Aerospace Engineering}, Vol. 2020, 2020, p.
  8813685.
\newblock \doi{10.1155/2020/8813685}.

\bibitem[{Chen et~al.(2015)Chen, Mahalec, Chen, Liu, He, and
  Sun}]{chen2015reconfiguration}
Chen, Y., Mahalec, V., Chen, Y., Liu, X., He, R., and Sun, K.,
  \enquote{Reconfiguration of satellite orbit for cooperative observation using
  variable-size multi-objective differential evolution,} \emph{European Journal
  of Operational Research}, Vol. 242, No.~1, 2015, pp. 10--20.
\newblock \doi{10.1016/j.ejor.2014.09.025}.

\bibitem[{de~Weck et~al.(2004)de~Weck, de~Neufville, and
  Chaize}]{deweck2004staged}
de~Weck, O.~L., de~Neufville, R., and Chaize, M., \enquote{Staged Deployment of
  Communications Satellite Constellations in Low Earth Orbit,} \emph{Journal of
  Aerospace Computing, Information, and Communication}, Vol.~1, No.~3, 2004,
  pp. 119--136.
\newblock \doi{10.2514/1.6346}.

\bibitem[{Arnas and Linares(2022)}]{arnas2022uniform}
Arnas, D., and Linares, R., \enquote{Uniform Satellite Constellation
  Reconfiguration,} \emph{Journal of Guidance, Control, and Dynamics}, Vol.~45,
  No.~7, 2022, pp. 1241--1254.
\newblock \doi{10.2514/1.G006514}.

\bibitem[{{Ferringer} et~al.(2009){Ferringer}, {Spencer}, and
  {Reed}}]{ferringer2009many}
{Ferringer}, M.~P., {Spencer}, D.~B., and {Reed}, P., \enquote{Many-objective
  reconfiguration of operational satellite constellations with the
  Large-Cluster Epsilon Non-dominated Sorting Genetic Algorithm-II,} \emph{2009
  IEEE Congress on Evolutionary Computation}, 2009, pp. 340--349.
\newblock \doi{10.1109/CEC.2009.4982967}.

\bibitem[{Davis(2010)}]{davis2010constellation}
Davis, J.~J., \enquote{Constellation reconfiguration: Tools and analysis,}
  Ph.D. thesis, Texas A\&M University, 2010.

\bibitem[{Fakoor et~al.(2016)Fakoor, Bakhtiari, and Soleymani}]{fakoor2016}
Fakoor, M., Bakhtiari, M., and Soleymani, M., \enquote{Optimal design of the
  satellite constellation arrangement reconfiguration process,} \emph{Advances
  in Space Research}, Vol.~58, No.~3, 2016, pp. 372--386.
\newblock \doi{10.1016/j.asr.2016.04.031}.

\bibitem[{Denis et~al.(2016)Denis, {de Boissezon}, Hosford, Pasco, Montfort,
  and Ranera}]{denis2016}
Denis, G., {de Boissezon}, H., Hosford, S., Pasco, X., Montfort, B., and
  Ranera, F., \enquote{The evolution of Earth Observation satellites in Europe
  and its impact on the performance of emergency response services,} \emph{Acta
  Astronautica}, Vol. 127, 2016, pp. 619--633.
\newblock \doi{10.1016/j.actaastro.2016.06.012}.

\bibitem[{Voigt et~al.(2016)Voigt, Giulio-Tonolo, Lyons, Kučera, Jones,
  Schneiderhan, Platzeck, Kaku, Hazarika, Czaran, Li, Pedersen, James, Proy,
  Muthike, Bequignon, and Guha-Sapir}]{voigt2016global}
Voigt, S., Giulio-Tonolo, F., Lyons, J., Kučera, J., Jones, B., Schneiderhan,
  T., Platzeck, G., Kaku, K., Hazarika, M.~K., Czaran, L., Li, S., Pedersen,
  W., James, G.~K., Proy, C., Muthike, D.~M., Bequignon, J., and Guha-Sapir,
  D., \enquote{Global trends in satellite-based emergency mapping,}
  \emph{Science}, Vol. 353, No. 6296, 2016, pp. 247--252.
\newblock \doi{10.1126/science.aad8728}.

\bibitem[{Wang et~al.(2021)Wang, Wu, Xing, and Pedrycz}]{wang2021agile}
Wang, X., Wu, G., Xing, L., and Pedrycz, W., \enquote{Agile Earth Observation
  Satellite Scheduling Over 20 Years: Formulations, Methods, and Future
  Directions,} \emph{IEEE Systems Journal}, Vol.~15, No.~3, 2021, pp.
  3881--3892.
\newblock \doi{10.1109/JSYST.2020.2997050}.

\bibitem[{Paek et~al.(2019)Paek, Kim, and de~Weck}]{paek2019optimization}
Paek, S.~W., Kim, S., and de~Weck, O., \enquote{Optimization of Reconfigurable
  Satellite Constellations Using Simulated Annealing and Genetic Algorithm,}
  \emph{Sensors}, Vol.~19, No.~4, 2019.
\newblock \doi{10.3390/s19040765}.

\bibitem[{McGrath and Macdonald(2019)}]{mcgrath2019general}
McGrath, C.~N., and Macdonald, M., \enquote{General Perturbation Method for
  Satellite Constellation Reconfiguration Using Low-Thrust Maneuvers,}
  \emph{Journal of Guidance, Control, and Dynamics}, Vol.~42, No.~8, 2019, pp.
  1676--1692.
\newblock \doi{10.2514/1.G003739}.

\bibitem[{Zhang et~al.(2021)Zhang, Zhang, Jiao, Baoyin, and Li}]{zhang2021}
Zhang, Z., Zhang, N., Jiao, Y., Baoyin, H., and Li, J., \enquote{Multi-Tree
  Search for Multi-Satellite Responsiveness Scheduling Considering Orbital
  Maneuvering,} \emph{IEEE Transactions on Aerospace and Electronic Systems},
  2021, pp. 1--1.
\newblock \doi{10.1109/TAES.2021.3129723}.

\bibitem[{Morgan et~al.(0)Morgan, McGrath, and de~Weck}]{morgan2023}
Morgan, S.~J., McGrath, C.~N., and de~Weck, O.~L., \enquote{Optimization of
  Multispacecraft Maneuvers for Mobile Target Tracking from Low Earth Orbit,}
  \emph{Journal of Spacecraft and Rockets}, Vol.~0, No.~0, 0, pp. 1--10.
\newblock \doi{10.2514/1.A35457},
  \urlprefix\url{https://doi.org/10.2514/1.A35457}.

\bibitem[{Appel et~al.(2014)Appel, Guelman, and Mishne}]{appel2014optimization}
Appel, L., Guelman, M., and Mishne, D., \enquote{Optimization of satellite
  constellation reconfiguration maneuvers,} \emph{Acta Astronautica}, Vol.~99,
  2014, pp. 166--174.
\newblock \doi{10.1016/j.actaastro.2014.02.016}.

\bibitem[{Legge~Jr(2014)}]{legge2014optimization}
Legge~Jr, R.~S., \enquote{Optimization and valuation of recongurable satellite
  constellations under uncertainty,} Ph.D. thesis, Massachusetts Institute of
  Technology, 2014.

\bibitem[{{de Weck} et~al.(2008){de Weck}, Scialom, and
  Siddiqi}]{deweck2008optimal}
{de Weck}, O.~L., Scialom, U., and Siddiqi, A., \enquote{Optimal
  reconfiguration of satellite constellations with the auction algorithm,}
  \emph{Acta Astronautica}, Vol.~62, No.~2, 2008, pp. 112--130.
\newblock \doi{10.1016/j.actaastro.2007.02.008}.

\bibitem[{L{\"u}ders(1961)}]{luders1961}
L{\"u}ders, R.~D., \enquote{Satellite networks for continuous zonal coverage,}
  \emph{ARS Journal}, Vol.~31, No.~2, 1961, pp. 179--184.
\newblock \doi{10.2514/8.5422}.

\bibitem[{L{\"u}ders and Ginsberg(1974)}]{luders1974}
L{\"u}ders, R., and Ginsberg, L., \enquote{Continuous zonal coverage-a
  generalized analysis,} \emph{Mechanics and Control of Flight Conference},
  1974, p. 842.
\newblock \doi{10.2514/6.1974-842}.

\bibitem[{Walker(1970)}]{walker1970}
Walker, J.~G., \enquote{Circular orbit patterns providing continuous whole
  earth coverage,} Tech. rep., Royal Aircraft Establishment Farnborough (United
  Kingdom), 1970.

\bibitem[{Walker(1977)}]{walker1977}
Walker, J.~G., \enquote{Continuous whole-earth coverage by circular-orbit
  satellite patterns,} Tech. rep., Royal Aircraft Establishment Farnborough
  (United Kingdom), 1977.

\bibitem[{Walker(1984)}]{walker1984}
Walker, J.~G., \enquote{Satellite constellations,} \emph{Journal of the British
  Interplanetary Society}, Vol.~37, 1984, pp. 559--572.

\bibitem[{Draim(1987)}]{draim1987common}
Draim, J.~E., \enquote{A common-period four-satellite continuous global
  coverage constellation,} \emph{Journal of Guidance, Control, and Dynamics},
  Vol.~10, No.~5, 1987, pp. 492--499.
\newblock \doi{10.2514/3.20244}.

\bibitem[{Lee et~al.(2020)Lee, Shimizu, Yoshikawa, and Ho}]{lee2020satellite}
Lee, H., Shimizu, S., Yoshikawa, S., and Ho, K., \enquote{Satellite
  Constellation Pattern Optimization for Complex Regional Coverage,}
  \emph{Journal of Spacecraft and Rockets}, Vol.~57, No.~6, 2020, pp.
  1309--1327.
\newblock \doi{10.2514/1.A34657}.

\bibitem[{Zhu et~al.(2010)Zhu, Li, and Baoyin}]{zhu2010satellite}
Zhu, K.-J., Li, J.-F., and Baoyin, H.-X., \enquote{Satellite scheduling
  considering maximum observation coverage time and minimum orbital transfer
  fuel cost,} \emph{Acta Astronautica}, Vol.~66, No.~1, 2010, pp. 220--229.
\newblock \doi{10.1016/j.actaastro.2009.05.029}.

\bibitem[{Mortari et~al.(2004)Mortari, Wilkins, and Bruccoleri}]{mortari2004}
Mortari, D., Wilkins, M.~P., and Bruccoleri, C., \enquote{The Flower
  Constellations,} \emph{The Journal of the Astronautical Sciences}, Vol.~52,
  No.~1, 2004, pp. 107--127.
\newblock \doi{10.1007/BF03546424}.

\bibitem[{Avenda{\~n}o et~al.(2013)Avenda{\~n}o, Davis, and
  Mortari}]{avendano2013}
Avenda{\~n}o, M.~E., Davis, J.~J., and Mortari, D., \enquote{The 2-D lattice
  theory of flower constellations,} \emph{Celestial Mechanics and Dynamical
  Astronomy}, Vol. 116, No.~4, 2013, pp. 325--337.
\newblock \doi{10.1007/s10569-013-9493-8}.

\bibitem[{Bartholdi et~al.(1980)Bartholdi, Orlin, and Ratliff}]{bartholdi1980}
Bartholdi, J.~J., Orlin, J.~B., and Ratliff, H.~D., \enquote{Cyclic Scheduling
  via Integer Programs with Circular Ones,} \emph{Operations Research},
  Vol.~28, No.~5, 1980, pp. 1074--1085.
\newblock \doi{10.1287/opre.28.5.1074}.

\bibitem[{Bartholdi(1981)}]{bartholdi1981}
Bartholdi, J.~J., \enquote{A Guaranteed-Accuracy Round-off Algorithm for Cyclic
  Scheduling and Set Covering,} \emph{Operations Research}, Vol.~29, No.~3,
  1981, pp. 501--510.
\newblock \doi{10.1287/opre.29.3.501}.

\bibitem[{Lee and Ho(2020)}]{Lee2020binary}
Lee, H., and Ho, K., \enquote{Binary Integer Linear Programming Formulation for
  Optimal Satellite Constellation Reconfiguration,} \emph{AAS/AIAA
  Astrodynamics Specialist Conference}, 2020.

\bibitem[{ReVelle et~al.(2008)ReVelle, Scholssberg, and Williams}]{revelle2008}
ReVelle, C., Scholssberg, M., and Williams, J., \enquote{Solving the maximal
  covering location problem with heuristic concentration,} \emph{Computers \&
  Operations Research}, Vol.~35, No.~2, 2008, pp. 427--435.
\newblock \doi{10.1016/j.cor.2006.03.007}, part Special Issue: Location
  Modeling Dedicated to the memory of Charles S. ReVelle.

\bibitem[{Megiddo et~al.(1983)Megiddo, Zemel, and Hakimi}]{megiddo1983}
Megiddo, N., Zemel, E., and Hakimi, S.~L., \enquote{The Maximum Coverage
  Location Problem,} \emph{SIAM Journal on Algebraic Discrete Methods}, Vol.~4,
  No.~2, 1983, pp. 253--261.
\newblock \doi{10.1137/0604028}.

\bibitem[{Church and {ReVelle}(1974)}]{church1974}
Church, R., and {ReVelle}, C., \enquote{The Maximal Covering Location Problem,}
  \emph{Papers in Regional Science}, Vol.~32, No.~1, 1974, pp. 101--118.
\newblock \doi{10.1111/j.1435-5597.1974.tb00902.x}.

\bibitem[{Avenda{\~{n}}o and Mortari(2009)}]{avendano2009}
Avenda{\~{n}}o, M., and Mortari, D., \enquote{{A closed-form solution to the
  minimum ${\Delta V_{\rm tot}^2}$ Lambert's problem},} \emph{Celestial
  Mechanics and Dynamical Astronomy}, Vol. 106, No.~1, 2009, p.~25.
\newblock \doi{10.1007/s10569-009-9238-x}.

\bibitem[{Vallado(2013)}]{vallado2013fundamentals}
Vallado, D., \emph{Fundamentals of Astrodynamics and Applications}, Space
  technology library, Microcosm Press, 2013, Chap.~6.

\bibitem[{Hoffman and Kruskal(1956)}]{hoffman1956}
Hoffman, A.~J., and Kruskal, J.~B., \emph{Integral Boundary Points of Convex
  Polyhedra}, Princeton University Press, 1956, Vol.~38, pp. 223--246.
\newblock \doi{10.1515/9781400881987-014}.

\bibitem[{Kuhn(1955)}]{kuhn1955}
Kuhn, H.~W., \enquote{The Hungarian method for the assignment problem,}
  \emph{Naval Research Logistics Quarterly}, Vol.~2, No. 1‐2, 1955, pp.
  83--97.
\newblock \doi{10.1002/nav.3800020109}.

\bibitem[{Bertsekas(1981)}]{Bertsekas1981}
Bertsekas, D.~P., \enquote{{A new algorithm for the assignment problem},}
  \emph{Mathematical Programming}, Vol.~21, No.~1, 1981, pp. 152--171.
\newblock \doi{10.1007/BF01584237}.

\bibitem[{Bertsekas and Eckstein(1988)}]{Bertsekas1988}
Bertsekas, D.~P., and Eckstein, J., \enquote{{Dual coordinate step methods for
  linear network flow problems},} \emph{Mathematical Programming}, Vol.~42,
  No.~1, 1988, pp. 203--243.
\newblock \doi{10.1007/BF01589405}.

\bibitem[{Yacov et~al.(1971)Yacov, Lasdon, and Wismer}]{haimes1971}
Yacov, H.~Y., Lasdon, L.~S., and Wismer, D.~A., \enquote{On a Bicriterion
  Formulation of the Problems of Integrated System Identification and System
  Optimization,} \emph{IEEE Transactions on Systems, Man, and Cybernetics},
  Vol.~1, No.~3, 1971, pp. 296--297.
\newblock \doi{10.1109/TSMC.1971.4308298}.

\bibitem[{Fisher(2004)}]{fisher2004}
Fisher, M.~L., \enquote{The Lagrangian Relaxation Method for Solving Integer
  Programming Problems,} \emph{Management Science}, Vol.~50, No.
  12\_supplement, 2004, pp. 1861--1871.
\newblock \doi{10.1287/mnsc.1040.0263}.

\bibitem[{Held and Karp(1971)}]{held1971}
Held, M., and Karp, R.~M., \enquote{{The traveling-salesman problem and minimum
  spanning trees: Part II},} \emph{Mathematical Programming}, Vol.~1, No.~1,
  1971, pp. 6--25.
\newblock \doi{10.1007/BF01584070}.

\bibitem[{Held et~al.(1974)Held, Wolfe, and Crowder}]{held1974}
Held, M., Wolfe, P., and Crowder, H.~P., \enquote{Validation of subgradient
  optimization,} \emph{Mathematical Programming}, Vol.~6, No.~1, 1974, pp.
  62--88.
\newblock \doi{10.1007/BF01580223}.

\bibitem[{Bertsimas and Tsitsiklis(1997)}]{bertsimas-LPbook}
Bertsimas, D., and Tsitsiklis, J., \emph{Introduction to linear optimization},
  Athena Scientific, 1997.

\bibitem[{Guignard(2013)}]{guignard2013}
Guignard, M., \emph{Lagrangian Relaxation}, Springer US, Boston, MA, 2013, pp.
  845--860.
\newblock \doi{10.1007/978-1-4419-1153-7_1168}.

\bibitem[{Geoffrion(1974)}]{geoffrion1974}
Geoffrion, A.~M., \emph{Lagrangean relaxation for integer programming},
  Springer Berlin Heidelberg, Berlin, Heidelberg, 1974, pp. 82--114.
\newblock \doi{10.1007/BFb0120690}.

\bibitem[{Golkar and {Lluch i Cruz}(2015)}]{golkar2015}
Golkar, A., and {Lluch i Cruz}, I., \enquote{The Federated Satellite Systems
  paradigm: Concept and business case evaluation,} \emph{Acta Astronautica},
  Vol. 111, 2015, pp. 230--248.
\newblock \doi{10.1016/j.actaastro.2015.02.009}.

\bibitem[{Lee and Ho(2021)}]{lee2021lagrangian}
Lee, H., and Ho, K., \enquote{A Lagrangian Relaxation-Based Heuristic Approach
  to Regional Constellation Reconfiguration Problem,} \emph{AAS/AIAA
  Astrodynamics Specialist Conference}, 2021.

\bibitem[{Galvão and ReVelle(1996)}]{galvao1996}
Galvão, R.~D., and ReVelle, C., \enquote{A Lagrangean heuristic for the
  maximal covering location problem,} \emph{European Journal of Operational
  Research}, Vol.~88, No.~1, 1996, pp. 114--123.
\newblock \doi{10.1016/0377-2217(94)00159-6}.

\end{thebibliography}

\end{document}